\documentclass[11pt,english]{article}

\usepackage[margin = 2.2cm]{geometry}

\usepackage{amsthm}
\usepackage{amsmath}
\usepackage{amssymb}
\usepackage{mathtools}
\usepackage{graphicx}
\usepackage[pdftex,colorlinks,backref=page,citecolor=blue,bookmarks=false]{hyperref}
\usepackage{color}
\usepackage{thm-restate}
\usepackage{enumitem}
\usepackage{subcaption}
\usepackage{xspace}
\usepackage{pifont}

\usepackage{caption}% http://ctan.org/pkg/caption
\usepackage[square,sort,comma,numbers]{natbib} 
\usepackage{setspace}
\usepackage{cleveref}
\theoremstyle{plain}

\newtheorem*{theorem*}{Theorem}
\newtheorem{theorem}{Theorem}[section]
\crefname{theorem}{Theorem}{Theorems}
\Crefname{theorem}{Theorem}{Theorems}

\newtheorem*{lemma*}{Lemma}
\newtheorem{lemma}[theorem]{Lemma}
\crefname{lemma}{Lemma}{Lemmas}
\Crefname{lemma}{Lemma}{Lemmas}

\newtheorem*{claim*}{Claim}

\crefname{claim}{Claim}{Claims}
\Crefname{claim}{Claim}{Claims}

\newtheorem*{innerclaim*}{Claim}

\crefname{claim}{Claim}{Claims}
\Crefname{claim}{Claim}{Claims}

\crefname{proposition}{Proposition}{Propositions}
\Crefname{proposition}{Proposition}{Propositions}

\crefname{corollary}{Corollary}{Corollaries}
\Crefname{corollary}{Corollary}{Corollaries}

\newtheorem{conjecture}[theorem]{Conjecture}
\crefname{conjecture}{Conjecture}{Conjectures}
\Crefname{conjecture}{Conjecture}{Conjectures}

\newtheorem{question}[theorem]{Question}
\crefname{question}{Question}{Questions}
\Crefname{question}{Question}{Questions}

\crefname{observation}{Observation}{Observations}
\Crefname{observation}{Observation}{Observations}

\crefname{example}{Example}{Examples}
\Crefname{example}{Example}{Examples}

\newtheorem{remark}[theorem]{Remark}
\crefname{remark}{Remark}{Remarks}
\Crefname{remark}{Remark}{Remarks}

\theoremstyle{definition}
\newtheorem{problem}[theorem]{Problem}
\crefname{problem}{Problem}{Problems}
\Crefname{problem}{Problem}{Problems}

\newtheorem{definition}[theorem]{Definition}
\crefname{definition}{Definition}{Definitions}
\Crefname{definition}{Definition}{Definitions}

\usepackage{xpatch}
\xpatchcmd{\proof}{\itshape}{\normalfont\proofnamefont}{}{}
\newcommand{\proofnamefont}{}
\renewcommand{\proofnamefont}{\bfseries}

\usepackage{framed}

\newcommand{\remove}[1]{}

\newcommand{\floor}[1]{
	\left\lfloor #1 \right\rfloor
}

\newcommand{\logstar}{\log^\star}

\newcommand{\cH}{\mathcal{H}}
\newcommand{\cM}{\mathcal{M}}

\newcommand{\cG}{\mathcal{G}}

\newcommand{\cP}{\mathcal{P}}

\newcommand{\cC}{\mathcal{C}}

\newcommand{\ksExpander}{KS-expander\xspace}
\newcommand{\ksExpanders}{KS-expanders\xspace}
\newcommand{\ssExpander}{SS-expander\xspace}
\newcommand{\ssExpanders}{SS-expanders\xspace}
\newcommand{\mExpander}{M-expander\xspace}
\newcommand{\mExpanders}{M-expanders\xspace}
\newcommand{\bmExpander}{BM-expander\xspace}
\newcommand{\bmExpanders}{BM-expanders\xspace}
\newcommand{\lExpander}{L-expander\xspace}
\newcommand{\lExpanders}{L-expanders\xspace}
\newcommand{\stExpander}{ST-expander\xspace}
\newcommand{\stExpanders}{ST-expanders\xspace}
\newcommand{\abExpander}{ABSZZ-expander\xspace}

\newcommand{\eps}{\varepsilon}
\newcommand{\lam}{\lambda}

\renewcommand{\setminus}{-}

\newcommand{\exr}{\ex^*}

\DeclareMathOperator{\diam}{diam}
\DeclareMathOperator{\ex}{ex}
\DeclareMathOperator{\sep}{sep}

\newcommand{\diK}{\overrightarrow{K}}
\newcommand{\TK}{\mathtt{TK}}
\renewcommand{\sb}{\mathrm{sb}}

\newcommand{\minor}{\mathrm{minor}}
\newcommand{\sub}{\mathrm{sub}}
\newcommand{\immers}{\mathrm{imm}}

\begin{document}

\title{Sublinear expanders and their applications}
\author{Shoham Letzter\thanks{
		Department of Mathematics, 
		University College London, 
		London WC1E~6BT, UK. 
		Email: \texttt{s.letzter}@\texttt{ucl.ac.uk}. 
		Research supported by the Royal Society.
	}
}

\date{}
\maketitle

\begin{abstract}

	\setlength{\parskip}{\medskipamount}
    \setlength{\parindent}{0pt}
    \noindent

	In this survey we aim to give a comprehensive overview of results using sublinear expanders. The term \emph{sublinear expanders} refers to a variety of definitions of expanders, which typically are defined to be graphs $G$ such that every not-too-small and not-too-large set of vertices $U$ has neighbourhood of size at least $\alpha |U|$, where $\alpha$ is a function of $n$ and $|U|$. This is in contrast with \emph{linear expanders}, where $\alpha$ is typically a constant.
	We will briefly describe proof ideas of some of the results mentioned here, as well as related open problems.

\end{abstract}

\section{Introduction} \label{sec:intro}

	Very informally speaking, expanders are graphs which have good connectivity properties, yet may be quite sparse. 
	Since their introduction by Bassalygo and Pinsker \cite{bassalygo1973complexity} in the 1970s, expanders have been studied extensively, and have seen numerous applications in combinatorics and computer science (see the concise expository paper by Sarnak \cite{sarnak2004expander} and the surveys by Hoory, Linial, and Wigderson \cite{hoory2006expander}, Lubotsky \cite{lubotzky2012expander}, and Krivelevich \cite{krivelevich2019expanders}).
	There is quite a large variety of definitions of expanders, but here we think of expanders as graphs whose every not-too-small and not-too-large set of vertices has a large neighbourhood. A very simple definition of expanders is the following, where $N_G(U)$, the \emph{neighbourhood} of $U$ in $G$, is the set of vertices in $G$ that are not in $U$ but have a neighbour in $U$. 

	\begin{definition}[Linear expanders] \label{def:linear-expander}
		Let $\alpha > 0$. A graph $G$ on $n$ vertices is called an \emph{$\alpha$-expander} if every subset $U \subseteq V(G)$ of size at most $n/2$ satisfies $|N_G(U)| \ge \alpha|U|$.
	\end{definition}

	When $\alpha$ is a constant, which is often the case in applications, this expansion property is \emph{linear}, namely, every not-too-large set of vertices expands linearly. See Krivelevich \cite{krivelevich2019expanders} for an excellent survey which mentions various results about $\alpha$-expanders and their variants, and illustrates various ways in which they can be used in applications.

	Our focus will be on \emph{sublinear expanders}. These may be defined similarly, but with $\alpha$ being a function of $n$ that tends to $0$ as $n$ grows, and oftentimes the rate of expansion of a set of vertices depends not only on $n$ but also on the size of the set. 
	Such expanders were first defined by Koml\'os and Szemer\'edi \cite{komlos1994topological,komlos1996topological} in 1994, who used them to solve a problem about subdivisions. A somewhat different notion of sublinear expanders was introduced by Shapira and Sudakov \cite{shapira2015small}, in their work about finding small minors. These two papers gave rise to a host of other results, particularly in the last three years or so.

	In this survey we aim to give a comprehensive overview of results that were proved using sublinear expanders, grouped into sections according the topics they address. We will also briefly describe proof ideas of some of the results presented here. See the survey of Liu \cite{liu2020robust} for a deeper dive into some of the methods developed for tackling problems using Koml\'os and Szemer\'edi's sublinear expanders.
	We will also highlight some open problems related to the results we will present.
	\begin{itemize}
		\item
			In \Cref{sec:subdivisions} we mention several results about average degree conditions implying the existence of a subdivision of a certain graph. 
		\item
			In \Cref{sec:small-minors-subdivisions} we mention a similar type of problem, about average degree conditions implying a small minor or subdivision of a complete graph. 
		\item
			\Cref{sec:immersion} discusses some results about immersions of complete graphs and digraphs. Next, we consider recent progress about the `odd cycle problem' of Erd\H{o}s and Hajnal, and a related problem about `balanced subdivisions' in \Cref{sec:odd-cycle-balanced-subdivision}.
		\item
			\Cref{sec:tight-cycle-rainbow-cycle-many-chords} presents several extremal results, about tight cycles, rainbow cycles and clique subdivisions, and cycles with many chords. 
		\item
			In \Cref{sec:decompositions-separation} we consider two problems about the global structure of a graph, namely about decomposing a graph into cycles and edges, and about separating the edges of a graph by paths.
		\item
			\Cref{sec:hamiltonian} mentions two results about the number of `Hamiltonian sets' in a graph with given average degree. 
		\item
			Finally, in \Cref{sec:other} we briefly mention a few other results that did not naturally belong in one of the previous sections.
	\end{itemize}

	\subsection*{Notation}

		Throughout the paper, we use $\log$ to denote the base 2 logarithm, and $\ln$ to denote the base $e$ logarithm.
		The various results mentioned in this survey use various different definitions of sublinear expanders. In order to easily distinguish the names of different expanders, we use the initials of those who introduced a specific definition of expanders; for example, the name \ksExpanders refers to expanders as defined by Koml\'os and Szemer\'edi.

		%\shoham{add numbers of lemmas/theorems in other papers}

\section{Subdivisions} \label{sec:subdivisions}
\begin{figure}[ht]
		\centering
		\includegraphics[scale = .6]{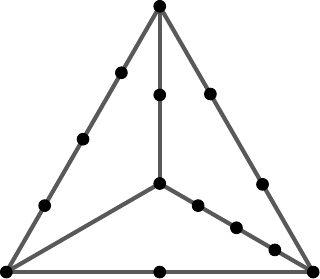}
		\caption{A subdivision of $K_4$}
		\label{fig:subdivision}
	\end{figure}

	Recall that, for a graph $F$, an \emph{$F$-subdivision} is a graph obtained by replacing each edge $uv$ in $F$ by a path with ends $u$ and $v$, such that the interiors of these paths are pairwise vertex-disjoint and vertex-disjoint of the original vertices of $F$ (see \Cref{fig:subdivision}).  

	Sublinear expanders were first introduced by Koml\'os and Szemer\'edi \cite{komlos1994topological} in 1994, who used them to solve an extremal problem about clique subdivisions. In this section we describe this first use, as well as several subsequent applications of sublinear expanders to similar problems about subdivisions. 

	\subsection{Finding clique subdivisions}
		
		In the 1990s Koml\'os and Szemer\'edi \cite{komlos1994topological,komlos1996topological}, as well as independently Bollob\'as and Thomason \cite{bollobas1996highly}, proved the following result, estimating the average degree that guarantees the existence of a $K_k$-subdivision.
		\begin{theorem}[Koml\'os--Szemer\'edi \cite{komlos1994topological,komlos1996topological} and Bollob\'as--Thomason \cite{bollobas1996highly}] \label{thm:Kk-subdivision}
			There exists $c > 0$ such that, for every integer $k \ge 1$, every graph with average degree at least $ck^2$ contains a $K_k$-subdivision.
		\end{theorem}
		This is tight up to the value of the constant $c$, as can be seen by considering a balanced complete bipartite graph on fewer than $2\binom{k/2}{2}$ vertices.

		The following notion of expanders played a key part in Koml\'os and Szemer\'edi's proof.\footnote{Bollob\'as and Thomason's proof was quite different and influential in its own right: their main result was that every $22k$-connected graph contains a graph which is \emph{$k$-linked}, namely for every sequence of distinct vertices $s_1,t_1, \ldots, s_k,t_k$, there is a collection of vertex-disjoint paths $P_1,\ldots,P_k$ such that $P_i$ joins $s_i$ with $t_i$.}
		\begin{definition}[Koml\'os--Szemer\'edi \cite{komlos1994topological,komlos1996topological}] \label{def:komlos-szemeredi}
			For $\eps, t > 0$, let $\rho(x) = \rho_{\eps,t}(x)$ be the function defined (for $x \ge t/2$).
			\begin{equation*}
				\rho(x) = \rho(x, \eps, t) = \frac{\eps}{(\log(15x/t))^2}.
			\end{equation*}
			An \emph{$(\eps, t)$-\ksExpander} is a graph $G$ in which every set of vertices $U$, with $t/2 \le |U| \le |G|/2$, satisfies $|N_G(U)| \ge \rho(|U|) \cdot |U|$.
		\end{definition}

		As this is a somewhat cumbersome definition, let us digest it briefly. Writing $n = |G|$, we remark that, typically, $t$ is much smaller than $n$; often we think of $t$ as constant and $n$ as large. Moreover, notice that every set of vertices of size $\Theta(t)$, and at least $t/2$, expands linearly, and sets of size $\Theta(n)$, and at most $n/2$, may expand at a rate as low as $O\big(\frac{1}{(\log n)^2}\big)$. Nevertheless, we have no information at all about the expansion rate of sets of fewer than $t/2$ vertices.

		At this point, it may not be at all clear why one may want to consider expanders as in \Cref{def:komlos-szemeredi}, when linear expanders, as in \Cref{def:linear-expander}, are much simpler to understand and work with. The basic reason is the fact that the sublinear $(\eps,t)$-\ksExpanders may be found in essentially any graph.

		\begin{theorem}[Koml\'os--Szemer\'edi \cite{komlos1994topological,komlos1996topological}] \label{thm:existence-expanders}
			Let $\eps > 0$ be sufficiently small, and let $t > 0$. Then every graph $G$ has a subgraph $H$ which is an $(\eps, t)$-\ksExpander, and satisfies $d(H) \ge d(G)/2$ and $\delta(H) \ge d(H)/2$.
		\end{theorem}

		\begin{remark} \label{rem:ks-expanders}
			\hfill
			\begin{itemize}
				\item
					The bound $d(H) \ge d(G)/2$ could be replaced by $d(H) \ge (1 - \delta)d(G)$, for any constant $\delta > 0$ using a slight variation of the proof (see, e.g., \Cref{thm:existence-expanders-robust}).
				\item
					We note that the definition of $\rho(x)$ is somewhat arbitrary. For the proof of \Cref{thm:existence-expanders} from \cite{komlos1996topological} to work, the following conditions need to hold: $\rho(x)$ is decreasing for $x \ge t/2$; $x\rho(x)$ is increasing for $x \ge t/2$; and $\int_{t/2}^{\infty} \frac{\rho(x)}{x}dx$ is finite. The exact choice of $\rho(x)$ was chosen for ease of presentation: the largest $\rho(x)$ can be while satisfying these conditions is $\rho(x) = \frac{1}{\log x (\log\log x)^{\Theta(1)}}$, which is close to best possible, as shown by Moshkovitz and Shapira \cite{moshkovitz2018decomposing}.
				\item
					Additionally, notice that $|H|$ can be much smaller than $|G|$ (e.g.\ if $G$ is a disjoint union of $K_{d+1}$'s, then $|H| \le d+1$, while $|G|$ can be arbitrarily large); this is sometimes inconvenient in applications.
				\item	
					Finally, note that there is some freedom in the choice of $t$. A natural choice is to take $t$ to be a small constant times $d(G)$ (the average degree of $G$), because then the lower bound on the neighbourhood of any given vertex is large enough for this neighbourhood to be guaranteed to expand by \Cref{def:komlos-szemeredi}.
			\end{itemize}
		\end{remark}

		A useful property of $(\eps, t)$-expanders is that any two not-too-small sets of vertices can be joined by a relatively short path, namely of length $O\big((\log n)^3\big)$.

		\begin{theorem}[Small diameter result; Koml\'os--Szemer\'edi \cite{komlos1994topological,komlos1996topological}] \label{thm:small-diameter}
			Let $G$ be an $n$-vertex $(\eps, t)$-expander. Then for every $x \ge t/2$ and every three sets of vertices $U_1, U_2, W$, where $|U_1|, |U_2| \ge x$ and $|W| \le \rho(x)x/4$, there is a path in $G - W$ between $U_1$ and $U_2$ of length at most $\frac{2}{\eps}\big(\log(15n/t)\big)^3$.
		\end{theorem}

		To see how \ksExpanders can be used to prove \Cref{thm:Kk-subdivision}, notice first that by \Cref{thm:existence-expanders}, it suffices to show that every $n$-vertex $(\eps, t)$-\ksExpander $H$ with average degree at least $c_1k^2$ and with $t = c_2k^2$, for some constants $c_1,c_2$, contains a $K_k$-subdivision.

		To prove this, Koml\'os and Szemer\'edi proved it separately for dense graphs, namely when $t = \Theta(n)$, where $n = |H|$, using Szemer\'edi's regularity lemma \cite{szemeredi1975regular}. 
		If $H$ is not dense, write $\diam = \frac{2}{\eps}(\log(15n/t))^3$, and notice that, roughly speaking, \Cref{thm:small-diameter} shows that the diameter of $H$ is close to $\diam$.
		In this case, if there are at least $k$ vertices with degree at least $D$, where $\binom{k}{2} \diam \le \rho(D) \cdot D/4$, then it is very easy to find a $K_k$-subdivision: given a set $K$ of $k$ vertices with degree at least $D$, apply \Cref{thm:small-diameter} for each pair of vertices in $K$ (with $U_1$ and $U_2$ being the neighbourhoods of these vertices, and $W$ the set of vertices used in previous connections) to find a path of length at most $\diam$ joining these vertices and avoiding previously used vertices, one by one. 

		Otherwise, the authors first apply \Cref{thm:existence-expanders} to the graph obtained by removing vertices of degree at least $D$, to find a new expander $H'$ with maximum degree less than $D$. Then, they find a set $K$ of $2k$ vertices, and associate a set $S(x)$ (called a `stable neighbourhood') with each $x \in K$, which expands well around $x$, and moreover no vertex appears in too many sets $S(x)$ (this last point is where they use that $H'$ has bounded degree). Now they again repeatedly apply \Cref{thm:small-diameter} to join pairs of vertices in $K$ through their stable neighbourhoods $S(\cdot)$, with a path of length at most $\diam$ avoiding previously defined paths, but whenever a set $S(x)$ becomes overused, they give up on the vertex $x$. A counting argument shows that at most $k$ vertices are discarded, and so this yields a $K_k$-subdivision.

		This idea of using sets $S(x)$ that are allowed to overlap somewhat, and starting with a bit more than the required $k$ vertices, allowed the authors to prove the tight $O(k^2)$ bound in \cite{komlos1996topological}; a similar argument where the sets $S(x)$ were required to be disjoint was used in their earlier paper \cite{komlos1994topological}, where a slightly weaker bound was proved.
		This idea, of allowing the $S(x)$ to overlap and then omitting $x$ whenever $S(x)$ becomes overused will be used in many subsequent results, with $S(x)$ being replaced by a variety of structures.

		While the bound $ck^2$ from \Cref{thm:Kk-subdivision} is tight up to a constant factor, it is still very interesting to find an asymptotically tight bound. Define $\sub(k)$ to be the minimum $d$ such that every graph with average degree at least $d$ contains a $K_k$-subdivision. We have seen that $\sub(k) = \Theta(k^2)$. The best known bound on $\sub(k)$ are
		\begin{equation*}
			\big(1 + o(1)\big) \cdot \frac{9k^2}{64} \le
			\sub(k) \le \big(1 + o(1)\big) \cdot \frac{10k}{23}.
		\end{equation*}
		The lower bound is an observation of {\L}uczak, using bipartite random graphs, and the upper bound is due to K\"uhn and Osthus \cite{kuhn2006extremal}.

		\begin{problem} \label{qn:subdivision}
			Determine $\sub(k)$ asymptotically. Is $\sub(k) = \big(1+o(1)\big) \frac{9k^2}{64}$?
		\end{problem}

	\subsection{Clique subdivisions in $C_4$-free graphs}

		While the requirement that $d(G) \ge ck^2$ is tight, up to a constant factor, for guaranteeing the existence of a $K_k$-subdivision, Mader \cite{mader1999extremal} conjectured that the quadratic bound could be replaced by a linear one for $C_4$-free graphs. Namely, he conjectured that there is a constant $c > 0$ such that, if $d(G) \ge ck$ and $G$ is $C_4$-free, then $G$ has a $K_k$-subdivision.

		In an early application of Koml\'os and Szemer\'edi's expanders, from 2004, K\"uhn and Osthus \cite{kuhn2004large} proved a slightly weaker bound: they showed that there is a constant $c > 0$ such that if $d(G) \ge ck(\log k)^{12}$ and $G$ is $C_4$-free, then $G$ contains a $K_k$-subdivision. In the same paper, K\"uhn and Osthus also prove an analogous result for $K_{s,t}$-free graphs, for all $s,t \ge 2$.

		In 2015, Balogh, Liu, and Sharifzadeh \cite{balogh2015subdivisions} proved Mader's conjecture under the additional assumption that $G$ is $C_6$-free; namely, they showed that there is a constant $c > 0$ such that, if $d(G) \ge ck$ and $G$ is $\{C_4, C_6\}$-free, then $G$ has a $K_k$-subdivision. 

		Mader's conjecture was subsequently solved by Liu and Montgomery \cite{liu2017proof} in 2017.
		\begin{theorem}[Liu--Montgomery \cite{liu2017proof}] \label{thm:mader}
			There is a constant $c > 0$ such that, if $d(G) \ge ck$ and $G$ is $C_4$-free, then $G$ contains a $K_k$-subdivision.
		\end{theorem}
		They also proved an analogous result for $K_{s,t}$-free graphs.

		We describe some elements of their proof of \Cref{thm:mader} here, some of which are inspired by \cite{balogh2015subdivisions} and also a result of Montgomery \cite{montgomery2015logarithmically} about small minors, that will be mentioned below.

		Here, the authors apply the existence result \Cref{thm:existence-expanders} with $t = \Theta(k^2)$. This yields an $n$-vertex $(\eps, t)$-expander $H$ with minimum degree at least $ck/4$. 
		By $C_4$-freeness, the second neighbourhood of every vertex has size at least $ck/4 \cdot (ck/4 -1)$, so if $t$ is taken to be a small factor of this number, then the expansion property of $H$ guarantees that second neighbourhoods expand well.

		The proof now splits into the following three cases, which are quite typical of a proof using $\ksExpanders$.
		\begin{itemize}
			\item
				There are at least $2k$ vertices of large degree (at least $k \cdot m^{c_1}$, where $m = \log(n/t)$ and $c_1$ is a large constant; this bound is chosen so that a greedy algorithm using the small diameter result \Cref{thm:small-diameter} would work).
			\item
				$H$ has small maximum degree and is quite dense ($d(G) \ge (\log n)^{c_2}$ for some large constant $c_2$).
			\item
				$H$ has small maximum degree and is sparse ($d(G) \le (\log n)^{c_2}$).
		\end{itemize}
		The first case is resolved similarly to \cite{komlos1996topological}: given a set $K$ of $2k$ large degree vertices, join pairs of these vertices one by one, using the small diameter property from \Cref{thm:small-diameter}, and removing vertices from $K$ whose neighbourhoods become overused. A counting argument shows that at least $k$ vertices remain in $K$, yielding the desired $K_k$-subdivision.

\begin{figure}[ht]
			\centering
			\begin{subfigure}{.3\textwidth}
				\includegraphics[scale = 1]{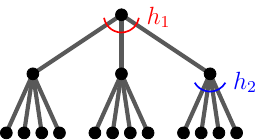}
			\end{subfigure}
			\hspace{.2cm}
			\begin{subfigure}{.6\textwidth}
				\centering
				\includegraphics[scale = .7]{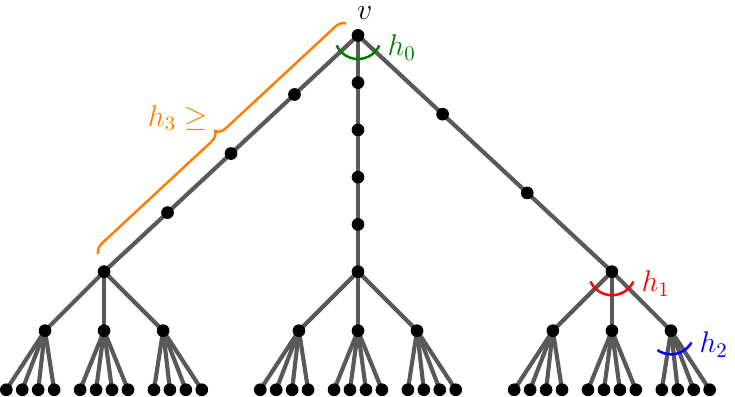}
			\end{subfigure}
			\caption{A $(3,4)$-hub and a $(5,3,4)$-unit}
			\label{fig:unit}
		\end{figure}

		If there are few vertices of large degree, they are removed from $H$, leaving a graph which is still an expander (with slightly worse parameters). 
		The second case is resolved somewhat similarly to the first one, except that instead of large degree vertices, structures called `units' are used; units are defined below, and inspired by a similar structure introduced by Montgomery \cite{montgomery2015logarithmically} (see \Cref{fig:unit}). The small maximum degree helps to find units that do not overlap too much.

		\begin{definition}[Hub] \label{def:hub}
			An \emph{$(h_1, h_2)$-hub} is a rooted tree of height $2$, where the root has degree $h_1$ and its neighbours have degree $h_2$.
		\end{definition}

		\begin{definition}[Unit] \label{def:unit}
			An \emph{$(h_0, h_1, h_2, h_3)$-unit} is a tree consisting of a core vertex $v$, $h_0$ pairwise vertex-disjoint $(h_1,h_2)$-hubs, and $h_0$ pairwise vertex-disjoint (except at $v$) paths of length at most $h_3$ joining the core vertex $v$ with each of the roots of the hubs.
		\end{definition}

		The authors of \cite{liu2017proof} take $h_0, h_1 = \Theta(k)$ and $h_2, h_3 = m^{\Theta(1)}$, where $m = \log(n/t)$, and then find $2k$ many $(h_0,h_1,h_2,h_3)$-units  whose non-leaf vertex sets are pairwise-disjoint. They construct them, one by one, by first constructing many disjoint hubs, and then joining their roots with short paths to find a unit. Notice that the requirement that $d(G)$ is relatively large is needed to guarantee the existence of many $(h_1,h_2)$-hubs, and the requirement that the maximum degree is not too large is used to show that the graph obtained by removing a not-too-large set of vertices still has large average degree, implying the existence of many hubs. With the $2k$ units at hand, one proceeds very similarly to the first case: connect, one by one, pairs of units, while avoiding vertices in any of the $2k$ units which are not leaves or their parents, and discard units whose sets of parents of leaves become overused.

		Finally, in the third case, using the sparsity and bounded maximum degree, they find $k$ vertices that are far from each other. They then connect pairs of vertices, one by one, by short paths that avoid previously chosen paths as well as the vicinity of other vertices.

	\subsection{Crux and clique subdivisions}

		Haslegrave, Hu, Kim, Liu, Luan, and Wang \cite{haslegrave2022crux} introduced a graph parameter, which they called `crux' and described as measuring the `essential order' of a graph, defined as follows.
		\begin{definition}[Crux] \label{def:crux}
			For a constant $\alpha \in (0,1)$ and a graph $G$, a subgraph $H \subseteq G$ is an \emph{$\alpha$-crux} if $d(H) \ge \alpha \cdot d(G)$. The $\alpha$-crux function of $G$, denoted $c_\alpha(G)$, is the minimum order of an $\alpha$-crux in $G$, namely,
			\begin{equation*}
				c_{\alpha}(G) = \min\{|H| : H \subseteq G \text{ and } d(H) \ge \alpha \cdot d(G)\}.
			\end{equation*}
		\end{definition}
		Thinking of $\alpha$ as a constant, notice that $c_{\alpha}(G) = \Omega(d(G))$. The authors of \cite{haslegrave2022crux} proved various generalisations of results about cycles in graphs with large average degree, showing that in many cases the average degree can be replaced by crux. For example, they showed that every graph $G$ contains a cycle of length $\Omega(c_{\alpha}(G))$, generalising (with a loss of a constant factor) Erd\H{o}s and Gallai's theorem \cite{erdos1959maximal} asserting that every graph $G$ contains a cycle of length at least $d(G)$.

		Liu and Montgomery \cite{liu2017proof} suggested that a common generalisation of \Cref{thm:mader} and \Cref{thm:Kk-subdivision}, that involves the average degree and crux\footnote{They did not use the term `crux' explicitly, as it was not introduced yet.} (yet yields slightly weaker bounds), might hold. Specifically, taking $\alpha = \frac{1}{100}$, say, they speculated that their methods could be used to show that every graph $G$ contains a $K_k$-subdivision, where
		\begin{equation*}
			k = \Omega\left(\!\min\left\{d(G), \sqrt{\frac{c_{\alpha}(G)}{\log c_{\alpha}(G)}}\,\right\}\right).
		\end{equation*}
		If true, this bound is tight: as noted in \cite{liu2017proof}, the $d$-blow-up $G$ of a $d$-vertex $O(1)$-regular expander satisfies $c_{\alpha}(G) = \Theta(d^2)$ and the largest clique subdivision has order $O\!\left(d \cdot (\log d)^{-1/2}\right)$.

		Im, Kim, Kim, and Liu \cite{im2022crux} proved a slightly weaker bound.
		\begin{theorem}[Im, Kim, Kim, and Liu \cite{im2022crux}] \label{thm:crux}
			Let $G$ be a graph, and write $t = \min\!\left\{d(G), \sqrt{\frac{c_{\alpha}(G)}{\log c_{\alpha}(G)}}\right\}$. Then $G$ contains a $K_k$-subdivision with $k = \Omega\big(t \cdot (\log \log t)^{-6}\big)$.
		\end{theorem}

		This paper is quite technical, and splits into several cases depending on density and value of $t$. In all cases the authors use vertex-disjoint stars, units whose non-leaf sets are disjoint, and webs, which are defined below, whose non-leaf sets are disjoint, to construct the desired subdivision. In some cases the crux is used to get an additional expansion property.
		
		It is natural to ask if $k$ can be taken to be $\Omega(t)$ in \Cref{thm:crux}.

		\begin{question}
			Is there a constant $c > 0$ such that for every graph $G$ if $t = \min\!\left\{d(G), \sqrt{\frac{c_{\alpha}(G)}{\log c_{\alpha}(G)}}\right\}$ then $G$ contains a $K_k$-subdivision with $k \ge c\cdot t$?
		\end{question}

	\subsection{Subdivisions of sparse graphs}

		So far we considered conditions guaranteeing a $K_t$-subdivision. Of course, one could instead look for $F$-subdivisions of other, sparser graphs $F$. In this direction, Haslegrave, Kim, and Liu \cite{haslegrave2022crux} considered $\alpha$-separable graphs: say that a graph $F$ is \emph{$\alpha$-separable} if there is a set $U$ of at most $\alpha|F|$ vertices, such that the components of $F \setminus U$ have size at most $\alpha |F|$. They proved the following for $\alpha$-separable graphs with bounded maximum degree.
		\begin{theorem}[Haslegrave--Kim--Liu \cite{haslegrave2022extremal}] \label{thm:separable}
			Let $\eps > 0$ and $\Delta \ge 1$, and let $\alpha$ be sufficiently small and $k$ sufficiently large. Then for every bipartite, $\alpha$-separable graph $F$ with $|F| \le (1 - \eps)k$, every graph $G$ with $d(G) \ge k$ contains an $F$-subdivision.
		\end{theorem}
		Notice that the assumption that $F$ is bipartite is crucial: if $F$ is a disjoint union of $1/\alpha$ complete graphs on $\alpha k / 2$ vertices, then $F$ is $\alpha$-separable, $|F| = k/2$, yet any complete bipartite graph $K_{t,t}$ with $t < \binom{\alpha k / 4}{2}$ does not contain a subdivision of $F$.
		The bound $|F| \le (1 - \eps)k$ is clearly optimal up to the error term $\eps k$ (consider $G = K_{k+1}$), but perhaps the error term could be decreased or even removed.
		This result has various implications regarding subdivisions of various sparse graphs, such as the grid, planar graphs, and minor-closed families.

		As part of their proof, the authors considered a `robust' notion of the expanders in \Cref{def:komlos-szemeredi}.
		\begin{definition}[Robust \ksExpanders] \label{def:komlos-szemeredi-robust}
			A graph $G$ is an \emph{$(\eps, t)$-robust-\ksExpander} if for every vertex set $U$ with $t/2 \le |U| \le |G|/2$, and every subgraph $F \subseteq G$ with $e(F) \le d(G) \cdot \rho(|U|) \cdot |U|$, the following holds.
			\begin{equation*}
				|N_{G \setminus F}(U)| \ge \rho(|U|) \cdot |U|.
			\end{equation*}
		\end{definition}

		They also proved an existence result, analogous to, and in fact slightly stronger than, \Cref{thm:existence-expanders}.
		(Here the notation $a \ll b$ means that there is a decreasing function $f$ such that the statement holds for $a \le f(b)$.)
		\begin{theorem} \label{thm:existence-expanders-robust}
			Let $0 < \nu \ll \eps_1, \eps_2 \ll \delta$.
			Then every graph $G$ with average degree $d$ has a subgraph $H$ which satisfies: $H$ is an $(\eps_1, \eps_2d)$-robust-\ksExpander; $d(H) \ge (1 - \delta)d$; $\delta(H) \ge d(H)/2$; and $H$ is $\nu d$-connected.
		\end{theorem}

\begin{figure}[ht]
			\centering
			\includegraphics[scale = .7]{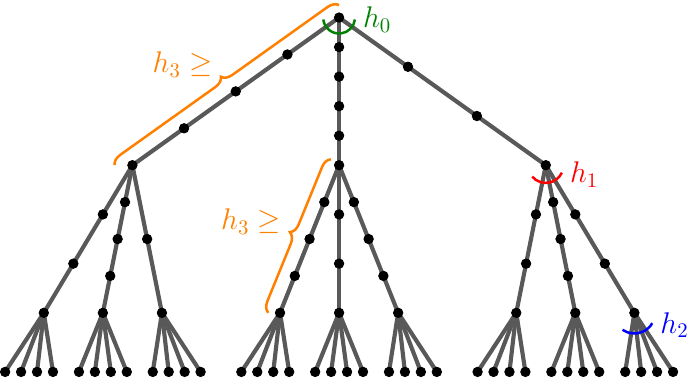}
			\caption{A $(3,3,4,5)$-web}
			\label{fig:web}
		\end{figure}

		This means, as usual, that it suffices to prove that every robust $(\eps, t)$-\ksExpander $H$, with $t = \Theta(k)$, $d(H) \ge k$, and $\delta(H) \ge k/2$, contains an $F$-subdivision, for $F$ as in \Cref{thm:separable}.
		And again, the proof differs depending on the density. The following three cases are considered.
		\begin{itemize}
			\item
				dense ($k = \Omega(n)$),
			\item
				medium ($k = o(n)$ and $k = \Omega\big( (\log n)^c\big)$ for some constant $c$),
			\item
				sparse ($k \le (\log n)^{O(1)}$).
		\end{itemize}
		In the dense case, the regularity lemma is used. In the medium case, when there are $k$ vertices of large enough degree (at least $k \cdot m^{c_1}$ for some constant $c_1$, where $m = \log(n/t)$), a greedy algorithm suffices. Otherwise, the authors show that the graph obtained by removing the large degree vertices still has large average degree, using the robust expansion of $H$ and the fact that $H$ is medium. They then proceed somewhat similarly to Liu and Montgomery \cite{liu2017proof}, finding `webs', which are reminiscent of `units' (see the next definition and \Cref{fig:web}).

		\begin{definition} \label{def:web}
			An \emph{$(h_0, h_1, h_2, h_3)$-web} is a rooted tree, obtained by subdividing the height 3 tree $T_{h_0, h_1, h_2}$ --- whose root has degree $h_0$, its neighbours have $h_1$ children, and these children in turn have $h_2$ children --- where the edges not touching leaves are subdivided at most $h_3$ times and edges touching leaves are not subdivided. The \emph{centre} of a web is the set of vertices on the paths which subdivide the edges touching the root in $T_{h_0,h_1,h_2}$.
		\end{definition}

		They then finds $2k$ webs, where $h_0, h_1, h_3$ are constant powers of $m = \log(n/t)$, and $h_0h_1h_2 \ge k \cdot m^{c_2}$ for some constant $c_2 > 0$, whose sets of non-leaf vertices are pairwise vertex-disjoint. They then connect them greedily, as usual, with connecting paths avoiding the centres of other webs, and omitting webs whose non-leaf sets become overused.

\begin{figure}[ht]
			\centering
			\includegraphics[scale = .8]{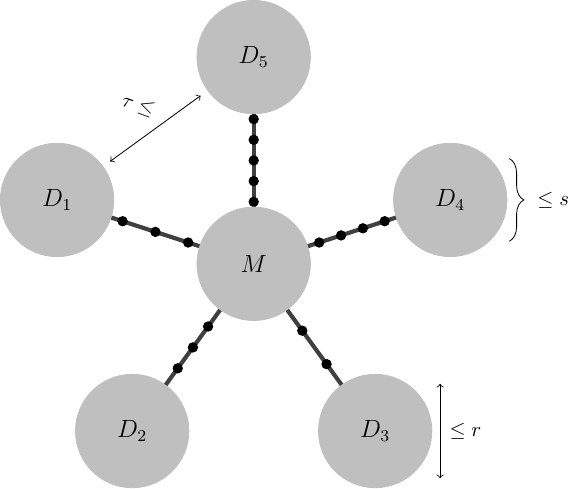}
			\caption{A nakji with $t = 5$. Its \emph{head} is $M$ and its \emph{legs} are $D_1, \ldots, D_5$}
			\label{fig:nakji}
		\end{figure}

		Finally, in the sparse case, the authors use `nakjis'\footnote{\emph{Nakji} means a `long arm octopus' in Korean.} (see \Cref{fig:nakji}).

		\begin{definition} \label{def:nakji}
			A \emph{$(t,s,r,\tau)$-nakji} in a graph $G$ is a subgraph $H$, consisting of vertex-disjoint sets $M$ and $D_i$, $i \in [t]$, and paths $P_i$, $i \in [t]$, such that the $P_i$'s are pairwise internally vertex-disjoint and their interiors are also disjoint of $M$ and the $D_j$'s; $P_i$ joins $M$ and $D_i$; $|D_i|, |M| \le s$; $D_i$ has diameter at most $r$; $M$ is $t$-connected; the sets $D_i$ and $M$ are at pairwise distance at least $\tau$. 

			We refer to $M$ as the \emph{head} of the nakji and to the $D_i$'s as the \emph{legs}\footnote{Perhas feet or hands would be more appropriate?}.
		\end{definition}

		In this case, the authors show that the graph obtained by removing large degree vertices still has average degree about as large as the original graph (otherwise, there is a copy of $F$ as a subgraph). Now, in this remaining graph, if there is a large, almost regular subexpander, then in it one can find enough vertices that are far from each other (using sparsity and almost regularity), and join them in a greedy fashion. If, instead, there is a subexpander which is medium or dense, then the arguments from the previous cases can be invoked. Otherwise, there are many small sparse subexpanders that are far from each other. They use these to find $k$ vertex-disjoint nakjis (with $t = \Delta$, i.e.\ the number of legs of each nakji corresponds to the maximum degree of $F$), where the head and legs of each nakji are distinct such subexpanders, which are joined together using the connectivity of the underlying expander $H$. 
		The nakjis are then connected, vertex-disjointly and through their legs, using the expansion properties of $H$. 
		The subdivision of $H$ is then finalised using the connectivity of the heads of the nakji to find appropriate star subdivisions.
		
\section{Small minors and subdivisions} \label{sec:small-minors-subdivisions}
	Say that a graph $G$ is an \emph{$H$-minor} if $H$ can be obtained from $G$ by removing vertices and edges and contracting edges.
	Equivalently, $G$ is an $H$-minor if there is a collection of vertex-disjoint connected subgraphs $(G_v)_{v \in V(H)}$, one for each vertex in $H$, such that if $uv$ is an edge in $H$ then there is an edge between $G_u$ and $G_v$. 

	Recall that $\sub(k)$ is the minimum $d$ such that every graph with average degree at least $d$ contains a $K_k$-subdivision.
	Define, analogously, $\minor(k)$ to be the minimum $d$ such that every graph with average degree at least $d$ contains a $K_k$-minor.
	We have seen that $\sub(k) = \Theta(k^2)$ (see \Cref{thm:Kk-subdivision}). Noting that a $K_k$-subdivision is also a $K_k$-minor, this implies that $\minor(k) = O(k^2)$. In fact, $\minor(k)$ is quite a lot smaller, and, moreover, its value is known quite precisely: Thomason \cite{thomason2001extremal} proved that $\minor(k) = \big(\alpha + o(1)\big) \cdot k\sqrt{\ln k}$, for an explicit constant $\alpha$.

	Fiorini, Joret, Theis, and Wood \cite{fiorini2012small} asked for the minimum $d$ guaranteeing an $H$-minor on \emph{few} vertices. Specifically, they asked if an average degree of at least $\minor(k) + \eps$ in an $n$-vertex graph guarantees an $H$-minor on $O_{\eps}(\log n)$ vertices. The size estimate would be tight, as can be seen by considering random graphs, and the average degree condition is tight by definition of $\minor(k)$. 

	\subsection{Almost logarithmically small clique minors}
		Addressing this question, in 2015 Shapira and Sudakov \cite{shapira2015small} proved the following.
		\begin{theorem}[Shapira--Sudakov \cite{shapira2015small}] \label{thm:minor-small-ss}
			For every $\eps > 0$ and integer $k \ge 1$, there is a constant $\alpha = \alpha(\eps, k)$ such that every $n$-vertex graph with average degree at least $\minor(k) + \eps$ has a $K_k$-minor on at most $\alpha \cdot \log n \cdot \log \log n$ vertices.
		\end{theorem}

		In their paper, they introduced a new notion of sublinear expanders.

		\begin{definition}[Shapira--Sudakov \cite{shapira2015small}] \label{def:shapira-sudakov}
			An $n$-vertex graph $G$ is a $\delta$-\emph{\ssExpander} if for every integer $d$ with $0 \le d \le \log \log n - 1$ and subset $S \subseteq V(G)$ of size at most $n/2^{2^d}$, we have
			\begin{equation*}
				|N(S)| \ge \frac{\delta \cdot 2^d}{\log n \cdot (\log \log n)^2} \cdot |S|.
			\end{equation*}
		\end{definition}

		To better understand this expression, for a set $S$ take $d$ to be largest such that $|S| \le n/2^{2^d}$, so $2^d \approx \log(n/|S|)$, yielding that $|N(S)| \gtrsim \frac{\delta \log(n/|S|)}{\log n \, (\log \log n)^2} \cdot |S|$. This shows that sets of size at most $n^{1 - \Omega(1)}$ expand at rate $\Omega(\frac{1}{(\log \log n)^2})$, and sets of size $\Theta(n)$ expand at rate $\Omega(\frac{1}{\log n \, (\log \log n)^2})$.

		As with \ksExpanders, every graph contains a $\delta$-\ssExpander, as long as $\delta$ is sufficiently small.

		\begin{theorem}[Shapira--Sudakov \cite{shapira2015small}] \label{thm:existence-ss-expanders}
			Let $\delta > 0$ be sufficiently small. Then every graph $G$ contains a $\delta$-\ssExpander $H$ with $d(H) \ge (1 - \delta)d(G)$.
		\end{theorem}

		An important reason why \ksExpanders are not appropriate here is that the small diameter result for them (\Cref{thm:small-diameter}) only yields a diameter of $\Theta( (\log n)^3)$, which is way too large for \Cref{thm:minor-small-ss}. Indeed, \ssExpanders do better in this respect. 
		We state a special case of it informally (see Claim 3.3 in \cite{shapira2015small} for a precise statement): in a $\delta$-\ssExpander on $n$ vertices, for every two sets $U_1, U_2$ of at least $(\log n)^4$ vertices, and a set $W$ of at most $(\log n)^2$ vertices, there is a path between $U_1$ and $U_2$ that avoids $W$ and has length $O(\log n \cdot (\log \log n)^3)$. This illustrates why this expansion notion is more suitable for the particular problem where we are interested in minors that are \emph{small}, as the diameter here is quite a lot smaller.
		However, using just this notion of expansion, Shapira and Sudakov get a slightly weaker result than the one stated above, showing that average degree $c(k) + \eps$ implies a $K_k$-minor on $O(\log n \, (\log \log n)^3)$ vertices. 

		Indeed, this is essentially immediate if there are at least $k$ vertices with degree at least $(\log n)^4$, using the `small diameter' result above, and in fact a $K_k$-subdivision is found in this case. If there are few vertices of high degree, the authors show that there are $k$ disjoint structures they call `expanding balls', which are relatively large sets of small radius and good expansion properties.
		They connect them via the same diameter argument, and show how to use these connections to get a small $K_k$-minor.

		To get the slightly stronger result, with a bound of $O(\log n \log \log n)$, they use yet another notion of expansion, and proved a corresponding existence result for this latter notion. We mention the definition and existence result here, as several future papers use a special case of the existence result, or variants of it.
		\begin{definition} \label{def:ss-expander-var}
			An $m$-vertex graph $H$ is said to be a \emph{$(\delta, n)$-\ssExpander} if for every integer $d$ with $0 \le d \le \log \log m - 1$ and $S \subseteq V(H)$ with $|S| \le m / 2^{2^d}$ we have $|N(S)| \ge \frac{\delta 2^d}{\log n} \cdot |S|$.
		\end{definition}

		\begin{theorem} \label{thm:ss-expander-var-existence}
			For every sufficiently small $\delta > 0$, every graph $G$ contains a subgraph $H$ such that $d(H) \ge (1-2\delta)d(G)$ and $H$ is a $(\delta, n)$-\ssExpander.
		\end{theorem}

	\subsection{Logarithmically small clique minors and subdivisions}

		Shortly afterwards, Montgomery \cite{montgomery2015logarithmically} improved upon the above result, proving the following tight result. 
		\begin{theorem}[Montgomery \cite{montgomery2015logarithmically}] \label{thm:small-minors-montgomery}
			For every $\eps > 0$ and integer $k \ge 1$, there is a constant $\alpha = \alpha(\eps,k)$ such that every $n$-vertex graph $G$ with $d(G) \ge \minor(k) + \eps$ has a $K_k$-minor on at most $\alpha \cdot \log n$ vertices.
		\end{theorem}

		He uses yet another definition of expanders. 
		\begin{definition}
			An $m$-vertex graph is a \emph{$(\lam, \eta)$-\mExpander} if every vertex set $S$ of size at most $m^{1-\eta}$ satisfies $|N(S)| \ge \lam |S|$.
		\end{definition}
		In his application, $\eta$ is a small constant, and $\lam$ a function of $\eta$, $m$, and the number of vertices $n$ in a given graph. Notice that this definition does not imply a small diameter result: sets of not-too-small size can be shown to expand quickly to size $m^{1-\eta}$, but might not expand fast beyond this size, making it hard to connect two given fixed vertices. 

		As usual, Montgomery's first step is to prove an existence result, showing that for given $\delta$ and $\eta$ and appropriate $\lam$, every $n$-vertex graph $G$ with average degree $d$ contains a subgraph $H$ with average degree at least $d(1 - \delta)$, which is either small or a $(\lam, \eta)$-\mExpander. Applying this to $G$ with average degree at least $\minor(t) + \eps$, we get a graph $H$ with $d(H) \ge \minor(t)$ which is either small or an expander. If $H$ is small, then in particular $|H| = O(\log n)$, and then it suffices to find a $K_t$-minor in $H$ (with no size restrictions), which is possible by the definition of $\minor(t)$. So suppose that $H$ is a $(\lam, \eta)$-\mExpander. Similarly to \cite{shapira2015small}, Montgomery shows that such $H$ contains many disjoint sets of size at least $m^{1/4}$ and small radius. Using the quick expansion of not-too-small sets to size $m^{1-\eta}$, he shows that there is a vertex $v$ which is at distance $O(\log n)$ from at least a $m^{-\eta}$ fraction of the sets, thereby overcoming the lack of a `small diameter' result. He uses the vertex $v$ as a basis of one of the $k$ vertices in the $K_k$-minor, and repeats the same argument after appropriate cleaning to find the desired $K_k$-minor.

		In the same paper, Montgomery also proved an analogous result for subdivisions. 
		\begin{theorem}[Montgomery \cite{montgomery2015logarithmically}] \label{thm:small-subdivision}
			For every $\eps > 0$ and integer $k \ge 1$, there is a constant $\alpha = \alpha(\eps,k)$ such that every $n$-vertex graph $G$ with $d(G) \ge \sub(k) + \eps$ has a $K_k$-subdivision on at most $\alpha \cdot \log n$ vertices.
		\end{theorem}

		For this, Montgomery used $(\lam, \eta)$-\mExpanders, whose every vertex set $S$ of size at most $m^{1/3}$ also expands at a rate of $\Omega(\frac{1}{(\log \log |S|)^2})$. This additional property is necessary here, as to find subdivisions we intuitively need single vertices to expand well. The proof here is more involved and uses the notion of `units', which in this context are a collection of $t$ disjoint, relatively large sets of small radius, along with $t$ paths from these sets to a common (`corner') vertex $v$, which are vertex-disjoint other than at $v$. This is somewhat similar to the notion of units from \Cref{def:unit}, and in fact served as inspiration for this latter definition of units.

\section{Immersions} \label{sec:immersion}

\begin{figure}[ht]
		\centering
		\includegraphics[scale = .7]{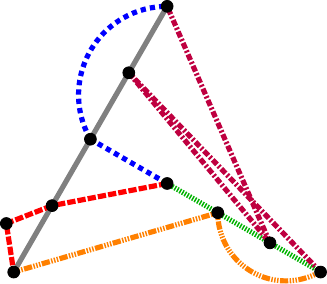}
		\caption{A $K_4$-immersion; each path $P_{uv}$ is depicted with a different pattern and colour.}
		\label{fig:immersion}
	\end{figure}

	Say that a graph $G$ \emph{immerses} a graph $H$ if there is an injective function $f : V(H) \to V(G)$, along with a collection of edge-disjoint paths $P_{uv}$, for $uv \in E(H)$, such that $P_{uv}$ has ends $u$ and $v$. Equivalently, $G$ immerses $H$ if $H$ can be obtained from $G$ by sequentially removing vertices and edges, and by replacing paths $uvw$ by the edge $uw$ (see \Cref{fig:immersion}).

	Recall that $\sub(k)$ is the minimum $d$ such that every graph with average degree at least $d$ contains a $K_k$-subdivision, and $\minor(k)$ is defined analogously for $K_k$-minors. We define $\immers(k)$ similarly, as the minimum $d$ such that every graph with average degree at least $d$ immerses $K_k$. Since a $K_k$-subdivision is also an immersion of $K_k$, we immediately get $\immers(k) \le \sub(k) = \Theta(k^2)$. Notice also that $\immers(k) > k-2$, because the existence of a $K_k$-immersion implies the existence of vertices of degree at least $k-1$.
	This was proved to be tight, up to a constant factor, by DeVos, Dvo\v{r}\'ak, Fox, McDonald, Mohar, and Scheide \cite{devos2014minimum}, who proved $\immers(k) \le 400k$. This was improved to $\immers(k) \le 22k + 14$ by Dvo\v{r}\'ak and Yepremyan \cite{dvorak2018complete}. The answer might actually be $k$, which would imply the following conjecture and answer the subsequent question affirmatively. While neither \cite{devos2014minimum} nor \cite{dvorak2018complete} use expanders, it is plausible that expanders can be used to make progress on the following conjecture and question.

	\begin{conjecture}[Lescure--Meyniel \cite{lescure1988problem}, Abu-Khzam--Langston \cite{abu2003graph}] \label{conj:immersion}
		If $\chi(G) \ge k$ then $G$ contains a $K_k$-immersion.
	\end{conjecture}

	\begin{question}[Dvo\v{r}\'ak--Yepremyan \cite{dvorak2018complete}] \label{qn:immersion}
		Does every graph with minimum degree $k$ immerse $K_k$?
	\end{question}

	We remark that Lescure and Meyniel \cite{lescure1988problem} and, independently, DeVos, Kawarabayashi, Mohar, and Okamura \cite{devos2010immersing} proved that minimum degree $k-1$ guarantees a $K_k$-immersion for $k \le 7$, thereby proving the above conjecture and answering the question affirmatively for this range of $k$'s, strengthening both conjecture and question slightly. Nevertheless, this is no longer true for large values of $k$ (see \cite{devos2014minimum,collins2014constructing}).

	Another natural question, given \Cref{thm:small-minors-montgomery} and \Cref{thm:small-subdivision} above, asks whether an average degree slightly above $\immers(k)$ suffices to guarantee a small immersion.

	\begin{question} \label{qn:small-immersion}
		Let $\eps > 0$ and $k \ge 1$.
		Is there a constant $\alpha = \alpha(\eps,k)$ such that every $n$-vertex graph $G$ with $d(G) \ge \immers(k)+\eps$ has a $K_k$-immersion on at most $\alpha \log n$ vertices?
	\end{question}

	\subsection{Clique immersions in $K_{s,t}$-free graphs}
		Liu, Wang, and Yang \cite{liu2022clique} proved \Cref{conj:immersion} asymptotically for $K_{s,t}$-free graphs, for every $s,t \ge 2$.

		\begin{theorem}[Liu--Wang--Yang \cite{liu2022clique}]
			For every $\eps > 0$, integers $s,t \ge 2$, and large enough $d$, every $K_{s,t}$-free graph $G$ with $d(G) \ge d$ contains a $K_k$-immersion with $k \ge (1 - \eps)d$.
		\end{theorem}
		
		The proof uses many of the tools described before. 
		First, by \Cref{thm:existence-expanders-robust}, it suffices to prove that every $n$-vertex $(\eps_1, \eps_2 k)$-\ksExpander $H$, which is $K_{s,t}$-free, and has average degree at least $(1 + \delta)k$, contains a $K_k$-immersion.
		Notice that $k$ cannot be too close to $n$, by the $K_{s,t}$-freeness. Thus, there are two main cases. 
		\begin{itemize}
			\item
				dense ($k \ge (\log n)^c$),
			\item
				sparse (otherwise).
		\end{itemize}
		In the dense case, the authors used a variant of `hubs' and `units' defined above: here an \emph{$(h_1, h_2, h_3)$-unit} consists of $h_1$ vertex-disjoint stars of size $h_2$, that are joined by edge-disjoint paths of length at most $h_3$ to a `core vertex'. Here edge-disjointness suffices because we are interested in immersions. They find $k$ edge-disjoint $(h_1, h_2, h_3)$-units (with $h_1 = k$, $h_2 = \diam^c$, and $h_3 = 2\diam$, where $\diam$ is the expression from the small diameter result \Cref{thm:small-diameter}) with the centres of stars pairwise disjoint, and then connect them greedily.

		The sparse case is somewhat more complex. If the maximum degree is at most $k \cdot (\log n)^c$, then one can find $k$ vertices with degree a bit over $k$ that are sufficiently far apart, and then join them, one by one, avoiding previously used edges and the vicinity of vertices not being currently joined. Also, as usual, if there are $k$ vertices of degree at least $k (\log n)^c$ then a greedy strategy can join them edge-disjointly to form an immersion. Thus, we may assume that there are fewer than $k$ vertices of degree at least $k(\log n)^c$. By $K_{s,t}$-freeness, the graph obtained by removing them still has large average degree. If the latter graph has a subexpander with average degree at bit above $d$ and which is either dense or has small maximum degree, then previous arguments can be applied. If not, this means that there are many subexpanders with average degree a bit above $k$, that are far from each other and each have a large degree vertex. Take $k$ such subexpanders and a largest degree vertex from each, and join these up as usual.

	\subsection{Immersions in directed graphs}

		Notice that the definition of immersions can be carried through to directed graphs.
		However, an analogue of $\immers(k)$ for digraphs does not exist: there are digraphs with arbitrarily large minimum in- and out-degree which do not immerse $\diK_3$ (the complete digraph on $3$ vertices); see Lochet \cite{lochet2019immersion}. Nevertheless, in the same paper Lochet showed that minimum out-degree $k^3$ guarantees an immersion of a transitive tournament on $\Omega(k)$ vertices, and it is plausible that minimum out-degree $k$ already suffices for an immersion of a transitive tournament on $\Omega(k)$ vertices.

		\begin{question}\label{qn:immersion-transitive}
			Is there a constant $c > 0$ such that, for every integer $k \ge 1$, every digraph with minimum out-degree at least $ck$ contains an immersion of the transitive tournament on $k$ vertices?
		\end{question}

		DeVos, McDonald, Mohar, and Scheide \cite{devos2012immersing,devos2013note} showed that every Eulerian digraph with minimum in-degree $k^2$ immerses $\diK_k$, and asked whether a linear bound would suffice.
		In \cite{girao2023immersion}, Gir\~ao and Letzter answered this question affirmatively. 
		\begin{theorem}[Gir\~ao--Letzter \cite{girao2023immersion}] \label{thm:immersion-eurlerian}
			There exists $c > 0$ such that every Eulerian digraph with minimum in-degree at least $ck$ immerses $\diK_k$.
		\end{theorem}

		The proof is a rare application of expanders in digraphs. We define an analogue of robust $(\eps, t)$-\ksExpanders (see \Cref{def:komlos-szemeredi-robust}), where instead of lower bounding the size of the neighbourhood $N(U)$, we lower bound the sizes of the out- and in-neighbourhoods of $U$. 
		To prove an existence result for such directed robust expanders, we use an undirected version guaranteeing the existence of an expander where every set of relevant size has large edge boundary (as opposed to the vertex boundary that was used above). Applying this to the underlying undirected graph obtained from the original Eulerian digraph $G$, which we additionally assume to be regular, we get a digraph $D' \subseteq G$ where every set of vertices $U$ of relevant size has many outgoing edges (i.e.\ edges from $U$ to $V(D') - U$) or many incoming edges (edges from $V(D') - U$ to $U$). 
		Now, we observe that because $G$ is Eulerian, it immerses an Eulerian \emph{multidigraph} $D$ that contains $D'$ as a subgraph.
		By the property of $D'$ mentioned above and by $D$ being Eulerian, every set of vertices $U$ in $D$ of relevant size has many outgoing edges \emph{and} many incoming edges, which readily implies that $D$ is a directed robust $(\eps,t)$-\ksExpander (using that $D$ is close to regular, which follows from regularity of $G$).

		Additionally, we use the fact that $\immers(k) = O(k)$, along with structural arguments, to show that every Eulerian multidigraph with minimum degree at least $ck$ and $O(k)$ vertices (with the technical condition that the underlying undirected graph has $\Omega(k^2)$ edges) immerses $\diK_k$.
		
		The last two paragraphs imply that it suffices to show that an $n$-vertex Eulerian expander $D$ (which in our setting is a multidigraph), with average degree at least $ck$, immerses a (simple) digraph on $\Theta(k)$ vertices with $\Omega(k^2)$ edges. Conveniently, due to the nature of immersions, we may assume that $D$ has maximum in- and out-degree $\Theta(k)$. As usual, we distinguish between dense ($k \ge (\log(n/k))^c$) and sparse expanders. In the sparse case we find $2k$ vertices which are far from each other, $k$ of which have large out-degree (not counting multiplicities) and the other $k$ having large in-degree. We then find edge-disjoint paths from each of the latter vertices to each of the former ones, thus forming an immersion of $\diK_{k,k}$. In the dense case, we follow a similar strategy to that of Koml\'os and Szemer\'edi \cite{komlos1996topological}, finding $2k$ large out-stars and $2k$ large in-stars, with distinct centres and where leaves are not shared by too many stars, and connect as many of the out-stars to the in-stars with short edge-disjoint directed paths as possible, yielding an immersion of a dense subgraph of $\diK_{2k,2k}$, as needed.

\section{The odd cycle problem and balanced subdivisions} \label{sec:odd-cycle-balanced-subdivision}

	In this section we mention recent developments about a conjecture of Erd\H{o}s and Hajnal about cycle lengths in graphs with large chromatic number, known as the `odd cycle problem', and about average degree conditions guaranteeing the existence of so-called balanced subdivisions.
	While seemingly unrelated at first glance, in both topics it is useful to be able to join two given vertices by a path of specific length. 

	\subsection{Cycle lengths in graphs with large average degree}
		For a graph $G$, let $\cC(G)$ be the set of cycle lengths in a graph $G$. In 1966, Erd\H{o}s and Hajnal \cite{erdos1966chromatic} suggested to study the sum $\sum_{\ell \in \cC(G)}\frac{1}{\ell}$ as a measure of the density of $G$'s cycle lengths. In particular, they asked whether $\sum_{\ell \in \cC(G)}\frac{1}{\ell}$ tends to infinity as $\chi(G) \to \infty$. In fact, Erd\H{o}s later \cite{erdos1975some} suggested that the same should hold as $d(G) \to \infty$, and, moreover, he thought it likely that $\sum_{\ell \in \cC(G)}\frac{1}{\ell} \ge \big(1/2 + o_d(1)\big)\ln  d$, for every graph $G$ with $d(G) \ge d$, which would be tight, as can be seen by considering $K_{d,d}$. In a slightly different direction, Erd\H{o}s \cite{erdos1984some} asked if every graph $G$ with sufficiently large average degree contains a cycle whose length is a power of $2$.

		In a breakthrough paper, Liu and Montgomery \cite{liu2023solution} answered both questions affirmatively. 
		\begin{theorem}[Liu--Montgomery \cite{liu2023solution}] \label{thm:odd-cycle}
			For large enough $d$, and every graph $G$ with $d(G) \ge d$ there exists $\ell \ge \frac{d}{10\, (\ln d)^{12}}$ such that $\cC(G)$ contains every even number in $[(\ln  \ell)^8, \ell]$.
		\end{theorem}

		They use similar methods to prove a similar result about the odd cycle lengths in a graph $G$ with large chromatic number, and, more generally, about cycle lengths of specific residues in graphs with large chromatic number. As an immediate corollary of their results, they solved a problem due to Erd\H{o}s and Hajnal, known as the `odd cycle problem': they asked whether $\cC_{\mathrm{odd}}(G)$, defined to be the sum of $1/\ell$ over all \emph{odd} cycle lengths of $G$, tends to infinity as $\chi(G) \to \infty$. Liu and Montgomery answer this affirmatively, and, in fact, they show that $\cC_{\mathrm{odd}}(G) \ge \big(1/2 + o(1)\big)\ln \chi(G)$.

\begin{figure}[ht]
			\centering
			\includegraphics[scale = .7]{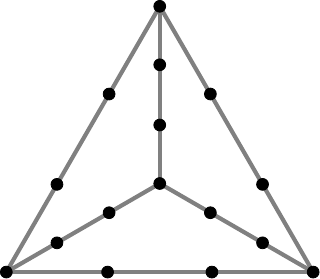}
			\caption{A balanced subdivision $\TK_4^{(3)}$ of $K_4$}
			\label{fig:bal-sub}
		\end{figure}

		In a somewhat different direction, they considered balanced subdivisions of complete graphs. A \emph{balanced subdivision} of a graph $H$ is a subdivision of $H$ where each edge is replaced by a path of the same length. Denote by $\TK_k^{(\ell)}$ the balanced subdivision of $K_k$ where each edge is replaced by a path of length $\ell$ (see \Cref{fig:bal-sub}).
		
		Liu and Montgomery proved that for any $k$, a large enough average degree guarantees the existence of a balanced $K_k$-subdivision, confirming a conjecture of Thomassen \cite{thomassen1984subdivisions}.
		\begin{theorem}[Liu--Montgomery \cite{liu2023solution}] \label{thm:balanced-subdivision-lm}
			For every integer $k \ge 1$ and large enough $d$, every graph $G$ with $d(G) \ge d$ contains a balanced subdivision of $K_k$.
		\end{theorem}

		While the results about cycle lengths and balanced subdivisions sound quite different, a common theme is the need to join pairs of vertices by a path of specific length. This is in contrast with previously described results, where when joining two given vertices, or sets of vertices, we just wanted the path joining them to be short, but did not care about the exact length. As such, the following theorem is a key component in the proofs of the results in \cite{liu2023solution}. Here an \emph{$(x,y)$-path} is a path with ends $x$ and $y$.

		\begin{theorem}[Theorem 2.7 in Liu--Montgomery \cite{liu2023solution}] \label{thm:diameter-exact}
			Let $\eps_1, \eps_2 > 0$ be small and let $d$ be large. Suppose that $H$ is a $\TK_{d/2}^{(2)}$-free bipartite $n$-vertex $(\eps_1, \eps_2 d)$-expander $H$ with $\delta(H) \ge d$. Then for any two vertices $x,y$ and $\ell \in [(\ln n)^7, n/(\ln n)^{12}]$ of the right parity\,\footnote{By the right parity, we mean that if $x,y$ are in the same part of the bipartition of $H$ then $\ell$ is even, and otherwise it is odd.}, there is an $(x,y)$-path of length $\ell$.
		\end{theorem}
		
		It is easy to see that this implies the first result about even cycle lengths. For the other results, some more work is needed.

		The main novelty in the proof of \Cref{thm:diameter-exact}, is the introduction and use of `adjusters'.
		\begin{definition}[Simple adjusters] \label{def:simple-adjusters}
			A \emph{simple $(D,m)$-adjuster} is a subgraph consisting of a cycle $C$ of length $2\ell$ for some $\ell \le 5m$, two vertices $v_1, v_2$ in $C$ at distance $\ell-1$ on $C$, and two pairwise vertex-disjoint graphs $F_1, F_2$, which are vertex-disjoint of $V(C) \setminus \{v_1, v_2\}$, such that, for $i \in [2]$: $v_i \in V(F_i)$; $|F_i| = D$; and every vertex in $F_i$ is at distance at most $m$ in $F_i$.
		\end{definition}
		\begin{figure}[ht]
			\centering
			\includegraphics[scale = .7]{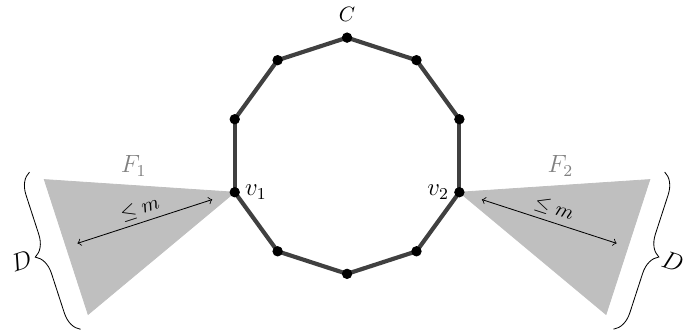}
			\caption{A simple adjuster}
			\label{fig:simple-adjuster}
		\end{figure}

		Assuming we have a long sequence of $k$ simple $(D,m)$-adjusters (where $m = \Theta((\log n)^3)$, so that $m$ is larger than the bound on the path length in \Cref{thm:small-diameter}, and $D$ is polylogarithmic in $n$ but significantly larger than $m$), the plan would be to join them up, one by one, via the sets $F_2$ and $F_1$ of consecutive adjusters, ensuring that the connecting paths are pairwise vertex-disjoint and relatively short. 
		This yields the following structure, generalising simple adjusters, which gives $(v_1,v_2)$-paths of many different lengths.
		\begin{definition}[Adjusters] \label{def:adjusters}
			A $(D,m,k)$-adjuster is a subgraph consisting of vertices $v_1,v_2$ and graphs $A,F_1,F_2$, such that: $A, F_1, F_2$ are pairwise vertex-disjoint; $v_i \in V(F_i)$ and every vertex in $F_i$ is at distance at most $m$ from $v_i$ in $F_i$; $|F_i| = D$ and $|A| \le 10mk$; and for some $\ell_0 \le 10mk$ and every $i \in \{0,1,\ldots,k\}$ there is a $(v_1,v_2)$-path in $A \cup \{v_1,v_2\}$ of length $\ell_0+2i$.
		\end{definition}
		\begin{figure}[ht]
			\centering
			\includegraphics[scale = .65]{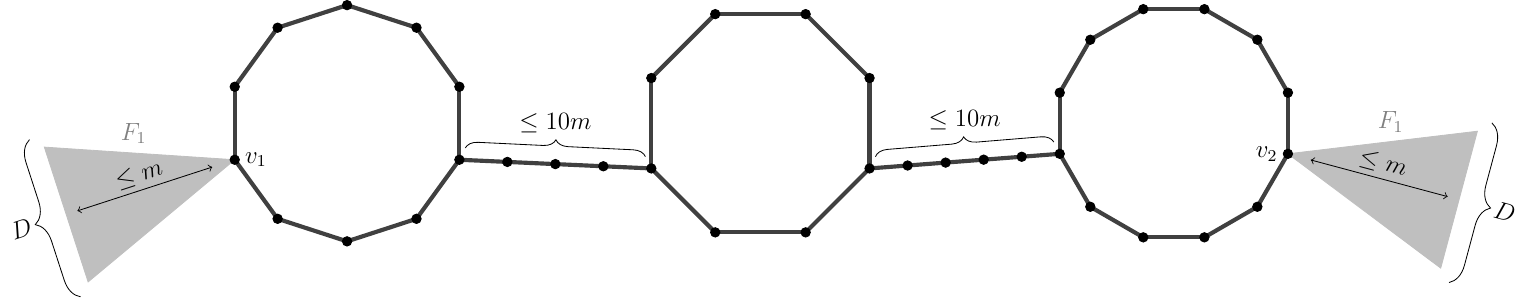}
			\caption{A $(D,m,3)$-adjuster}
			\label{fig:adjuster}
		\end{figure}

		Now, given two vertices $x,y$ and a $(D,m,k)$ adjuster, if we would like to find an $(x,y)$-path of length $\ell$, it suffices to join $x$ to $F_1$ and $y$ to $F_2$ by paths of total length between $\ell - \ell_0 - 2k$ and $\ell - \ell_0 - 2m$, which can then be corrected to an $(x,y)$-path of length $\ell$ using the adjuster. 
		These tasks, of joining the $F_i$'s to form a $(D,m,k)$-adjuster, and of joining $x$ and $y$ to the adjuster by a path of length approximately $\ell$, are both achievable via the diameter theorem \Cref{thm:small-diameter}, assuming that $k \ge m^c$ for large enough $c$. 

		As such, it suffices to show how to find many disjoint simple adjusters, which essentially amounts to finding a simple adjuster avoiding a forbidden set of vertices $W$ of polylogarithmic size. Finding one simple adjuster in an expander (with no forbidden vertices) is quite simple: take a shortest cycle $C$, pick two vertices on it at appropriate distance in $C$, and expand from there, showing that $C$ can be avoided due to it being a shortest cycle. To find an adjuster avoiding a set $W$, one can apply this reasoning to a subexpander of $G \setminus W$. However, as we have no control over the order of such a subexpander, this might lead to an adjuster which is way too small. Thus, the authors take many disjoint such adjusters, and show that one of them can be expanded to be sufficiently large. This is a challenging task, which we do not elaborate on. We do briefly mention that a useful tool in overcoming this challenge is a lemma (Lemma 3.7 in \cite{liu2023solution}) about expanding a set $A$ while avoiding forbidden sets with various ways of controlling the interaction between $A$ and the forbidden sets.

	\subsection{Improved bounds on average degree implying balanced clique subdivision}

		Notice that, through \Cref{thm:balanced-subdivision-lm}, Liu and Montgomery proved the existence of a function $\sb(k)$ such that every graph with average degree at least $\sb(k)$ contains a balanced subdivision of $K_k$. They do not calculate explicitly the upper bound on $\sb(k)$ that their proof yields, but they suggested that $\sb(k)$ might be $O(k^2)$. If true, this would be tight, and would generalise the aforementioned fact that $\sub(k) = \Theta(k^2)$ (see \Cref{thm:Kk-subdivision}).

		The first progress towards estimating $\sb(k)$ was made by Wang \cite{wang2021balanced}, who proved that $\sb(k) \le k^{2+o(1)}$. This was improved to the tight $\sb(k) = O(k^2)$ by Luan, Tang, Wang, and Yang \cite{luan2023balanced}, and, independently, by Gil Fernandez, Hyde, Liu, Pikhurko, and Wu \cite{gilfernandez2023disjoint}. 
		\begin{theorem}[Luan--Tang--Wang--Yang \cite{luan2023balanced} and Gil Fernandez--Hyde--Liu--Pikhurko--Wu \cite{gilfernandez2023disjoint}] 
			There is a constant $c > 0$ such that for every integer $k \ge 1$, every graph with average degree at least $ck^2$ contains a balanced subdivision of $K_k$.
		\end{theorem}
		The former paper \cite{wang2021balanced} also proves that there exists $c > 0$ such that average degree at least $ck$ suffices to guarantee a balanced $K_k$-subdivision for $C_4$-free graphs, generalising 
		\Cref{thm:mader}. 

		As usual, it suffices to prove that every $n$-vertex robust $(\eps_1, \eps_2 d)$-\ksExpander $H$ contains the required balanced $K_k$-subdivision, with $k = \Omega(\sqrt{d})$.
		Again as usual, both proofs split into three cases, depending on density: sparse ($d = (\log n)^{O(1)}$); medium ($d \ge (\log n)^{c}$ for some constant $c$ and $d = o(n)$); and dense ($d = \Omega(n)$).

		The dense case follows from a result of Alon, Krivelevich, and Sudakov \cite{alon2003turan}, that actually yields a $\TK_{k}^{(2)}$.
		The sparse case was already addressed by Wang \cite{wang2021balanced}, and Wang addressed it as follows. Note that we may assume $H$ is $\TK_k^{(2)}$-free. The proof uses tools due to Liu and Montgomery \cite{liu2023solution} about joining two vertices by a path of specific length in a $\TK_{k}^{(2)}$-free expander, and is otherwise quite routine. If there are enough vertices of large degree (at least $d(\log n)^c$ for some constant $c$) then they can be joined greedily. Otherwise, we may assume the maximum degree is small, allowing us to find many `core' vertices (in the same part of the bipartition) that are far enough from each other. These can then be joined, one by one, by paths of specific length, while avoiding previously used vertices and the vicinity of other core vertices. 

		It thus remains to address the medium case. 
		In \cite{luan2023balanced}, the authors find $\Theta(k)$ units (see \Cref{def:unit}) with at least $d^2m^c$ leaves (where $m = \log(n/d)$ and $c$ is some large constant), such that the sets of non-leaves are pairwise vertex-disjoint. They then proceed to join pairs of units, one by one, removing units that become overused. To yield a balanced subdivisions, the connection here needs to be of a specific length, and for this the authors use $(D,m,k)$-adjusters (see \Cref{def:adjusters}). Since here the adjusters need to be somewhat larger than in \cite{liu2023solution} ($D$ here has size $dm^{\Theta(1)}$, to be able to use them in connections while avoiding $O(d m^{O(1)})$ vertices, as opposed to $m^{\Theta(1)}$ in \cite{liu2017proof}), their construction is somewhat different here.
		The authors of \cite{gilfernandez2023disjoint}, follow a similar approach, using a variant of units, where a $(h_0,h_1,h_2)$-unit is defined to be a rooted tree of height $h_2+1$, such that the root has degree $h_0$, the vertices at level between $1$ and $h_2-1$ have exactly one child, and the vertices at level $h_2$ have $h_1$ children. Units are also used in place of the graphs $F_1,F_2$ in the definition of adjusters.

\section{Tight cycles, rainbow subdivisions, and cycles with many chords} \label{sec:tight-cycle-rainbow-cycle-many-chords}

	In this section we describe several Tur\'an type results, about tight cycles, rainbow subdivisions, and cycles with many chords. These results are connected to each other via the methods they use, which can be traced back to the paper of Shapira and Sudakov \cite{shapira2015small} about small minors. In particular, various notions of expansions are used here (though for the most part we will not give the exact definitions), but all are somewhat similar to the one in \Cref{def:ss-expander-var}.

	\subsection{Hypergraphs with no tight cycles}

\begin{figure}[ht]
			\centering
			\includegraphics[scale = .6]{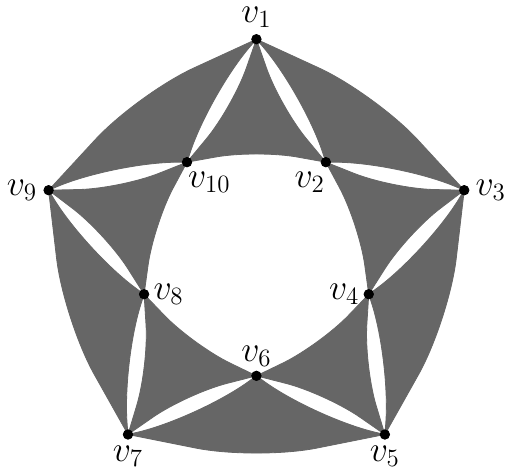}
			\caption{A $3$-uniform tight cycle on 10 vertices}
			\label{fig:tight-cycle}
		\end{figure}

		An $r$-uniform \emph{tight cycle} is a hypergraph on vertices $\{v_1, \ldots, v_{\ell}\}$, for some $\ell \ge r+1$, with edges $\{v_i \ldots v_{i+r-1}: i \in [\ell]\}$ (addition of indices taken modulo $r$; see \Cref{fig:tight-cycle}).
		
		It is natural to ask: what is the maximum number of edges an $n$-vertex $r$-uniform hypergraph can have, if it has no tight cycles? Denoting this maximum by $\ex_r(n, \cC)$, notice that $\ex_r(n, \cC) \ge \binom{n-1}{r-1}$ (take all edges containing a single vertex), which is best possible for $r = 2$. For larger values of $r$, this turns out not to be tight (disproving a conjecture due to S\'os and, independently, Verstra\"ete; see \cite{verstraete2016extremal,mubayi2011hypergraph}). The best known lower bound is due to Janzer \cite{janzer2020large}, who showed that $\ex_r(n, \cC) = \Omega\big(n^{r-1} \cdot \frac{\log n}{\log \log n}\big)$.
		Sudakov and Tomon \cite{sudakov2022extremal} proved a close-to-matching upper bound.
		\begin{theorem}[Sudakov--Tomon \cite{sudakov2022extremal}] \label{thm:sudakov-tomon}
			Let $r \ge 3$. Every $n$-vertex $r$-uniform hypergraph with no tight cycles has at most $n^{r-1}e^{O(\sqrt{\log n})}$ edges.
		\end{theorem}
		In particular, this shows $\ex_r(n, \cC) = n^{r-1+o(1)}$.
		
		They define $r$-line-graphs, which are graphs whose vertices are the edges of an $r$-partite $r$-uniform hypergraph, and whose edges correspond to two edges intersecting in $r-1$ vertices, and define an appropriate notion of density for such graphs. 
		They use a variant of a notion of expanders used by Shapira and Sudakov \cite{shapira2015small} (see \Cref{def:ss-expander-var}). 
		\begin{definition}[Sudakov--Tomon \cite{sudakov2022extremal}] \label{def:st-expanders}
			A \emph{$(\lam, d)$-\stExpander} is a graph $H$ which has minimum degree\footnote{Their definition of the minimum degree of an $r$-line-graph is somewhat different to the usual definition of the minimum degree. Denote by $\cH$ the underlying $r$-partite $r$-graph and by $(U_1, \ldots, U_r)$ an appropriate partition of $V(\cH)$. The \emph{minimum degree} of $H$ is defined to be the minimum, over all $i \in [r]$ and vertices $x$ in $H$, of the number of neighbours of $x$ in $H$ that agree with $x$ on $U_j$ for all $j \in [r] - \{i\}$ (plus $1$).} at least $d$, and where every vertex set $U \subseteq V(H)$ satisfies $|N(U)| \ge \lam|U|$.
		\end{definition}

		They then prove an existence result, showing that every $r$-line-graph $G$ of density $d$ (with an appropriate definition of density) contains a subgraph $H$ which is a $(\lam, d')$-\stExpander, with $\lam = \Omega(1/\log n)$ and $d' = \Omega(d)$, and has density at least $d/2$. This proof is inspired by similar arguments made by Shapira and Sudakov \cite{shapira2015small} to prove an existence result for expanders as in \Cref{def:ss-expander-var}. 

		Notice that for this notion of expansion to be effective in expanding a single vertex, one would like to have $d = \Omega(\log n)$. This is indeed the case here; in fact, we have $d = e^{\Omega(\sqrt{\log n})}$. A key property of expanders in $r$-line-graphs that is proved here is that for every vertex $e$ in an $(\lam,d)$-\stExpander $H$ (with appropriate parameters), which is an $r$-line-graph, one can reach almost every other vertex $f$ via a short tight path in the hypergraph corresponding to $H$. Observe that, because of the nature of tight paths, this does not imply that every two vertices can be joined via a short tight path, analogously to \Cref{thm:small-diameter}, even if we insist that these vertices correspond to disjoint edges.

		In order to find a tight cycle, the authors split the underlying vertex set of $H$ into two sets, denoting the two corresponding $r$-line-graphs by $H_1, H_2$. If they could prove that each graph $H_i$ maintains the expansion properties of $H$, then they would be able to find two tight paths, one from an $r$-set $e$ to another $r$-set $f$ in $H_1$, and one from $f$ to $e$ in $H_2$, using disjoint vertex sets (other than the vertices in $e,f$), thereby yielding a tight cycle, as required.
		The authors are unable to do so. Instead, they use the existence result from above to almost decompose each of $H_1$ and $H_2$ into expanders. If any of them are small, then we get a much denser expander than the original one, and the same argument is then repeated in this denser expander.
		Otherwise, the total number of expanders used in the decompositions is small, and then the expansion property can be used to essentially realise the strategy outlined above to get a tight cycle.
		Since the density cannot be increased indefinitely, at some iteration the latter option holds, yielding a tight cycle, as required.

		In \cite{letzter2023hypergraphs}, the author improved Sudakov and Tomon's bound, proving the almost tight bound $\ex_r(n, \cC) = O\big(n^{r-1}(\log n)^5\big)$. 
		\begin{theorem}[Letzter \cite{letzter2023hypergraphs}] \label{thm:hypergraph}
			Let $r \ge 3$. Every $n$-vertex $r$-uniform hypergraph with no tight cycles has at most $O(n^{r-1}(\log n)^5)$ edges.
		\end{theorem}

		Using the setup from \cite{sudakov2022extremal}, the improvement came from proving that every $r$-line-graph $H$, which is an \stExpander with appropriate parameters, and has sufficiently large density, has a tight cycle (so there is no need for a density increment argument as above). An important idea was to strengthen the expansion property above, showing that for every vertex $e$ in an $r$-line-graph $H$ as above, almost every vertex $f$ in $H$ can be reached from $e$ via a short tight path, in such a way that no vertex of the underlying graph is used in too many of these paths (except for the vertices in $e$). Another trick is then needed to overcome the directed nature of the problem.

		It would be interesting to prove an even tighter bound on $\ex_r(n,\cC)$. Recall that the best known lower bound \cite{janzer2020large} is $\ex_r(n,\cC) = \Omega\big(n^{r-1} \cdot \frac{\log n}{\log \log n}\big)$.
		\begin{problem}
			Determine $\ex_r(n,\cC)$ asymptotically. Is $\ex_r(n,\cC) = O\big(n^{r-1} \cdot \frac{\log n}{\log \log n}\big)$? A bit more crudely, is $\ex_r(n,\cC) = n^{r-1} \cdot (\log n)^{1 + o(1)}$?
		\end{problem}
		It is plausible that ideas described in subsequent subsections could yield a somewhat better than the one in \Cref{thm:hypergraph}, perhaps decreasing the exponent of $\log n$ (though likely they would require significant technical work), but it seems that to get a tight bound, new ideas will be needed.

	\subsection{Rainbow clique subdivisions}

		A similar notion of expansion was used by Jiang, Methuku, and Yepremyan \cite{jiang2023rainbow} to tackle a rainbow Tur\'an problem. A \emph{rainbow} graph is an edge-coloured graph whose edges have distinct colours. Keevash, Mubayi, Sudakov, and Verstra\"ete \cite{keevash2007rainbow} defined the \emph{rainbow Tur\'an number} of a graph $H$, denoted $\exr(n, H)$, as the maximum number of edges in an $n$-vertex properly edge-coloured graph which has no rainbow copies of $H$; a similar definition could be made for a collection of graphs $\cH$. 

		A particularly interesting question here is to determine the rainbow Tur\'an number of cycles, namely the maximum number of edges in an $n$-vertex properly edge-coloured graph which has no rainbow cycles, which we denote by $\exr(n, \cC)$. It is easy to see that the hypercube $Q_d$, where edges are coloured according to their direction, is a properly edge-coloured graph with no rainbow cycles, showing $\exr(n, \cC) = \Omega(n \log n)$. The authors of \cite{keevash2007rainbow} conjectured that the latter bound is tight, up to a constant factor.
		\begin{conjecture}[Keevash--Mubayi--Sudakov--Verstra\"ete \cite{keevash2007rainbow}]
			Every $n$-vertex properly edge-coloured graph with no rainbow cycles has $O(n \log n)$ edges.
		\end{conjecture}

		Janzer \cite{janzer2023rainbow} proved that $\exr(n, \cC) = O\big(n (\log n)^4\big)$, getting quite close to proving the conjecture. His main tool was an upper bound on the number of homomorphic copies of an even cycle with some degeneracy, and he used this method to prove several other Tur\'an type questions.

		The problem Jiang, Methuku, and Yepremyan \cite{jiang2023rainbow} considered was the following: for a constant $k$, what is the maximum number of edges in a properly edge-coloured $n$-vertex graph which has no rainbow subdivisions of $K_k$? Clearly, this number is at least $\exr(n, \cC) = \Omega(n \log n)$. The authors of \cite{jiang2023rainbow} proved the following bound, showing that the answer is $n^{1 + o(1)}$.
		\begin{theorem}
			Let $k \ge 3$. If $G$ is an $n$-vertex properly edge-coloured graph with no rainbow $K_k$-subdivision, then $e(G) \le n \cdot e^{O(\sqrt{\log n})}$.
		\end{theorem}
		An attentive reader might notice the similarity between this bound and that of Sudakov and Tomon \cite{sudakov2022extremal} in \Cref{thm:sudakov-tomon}. This is not an accident: the outline of the proof here is very similar to that in \cite{sudakov2022extremal}. The authors use a variant of expanders as in \Cref{def:st-expanders}, focusing on the edge boundary of vertex sets rather than the vertex boundary, and prove that in a properly-coloured expander, for every vertex $u$, almost every other vertex $v$ can be reached from $u$ via a short rainbow path. They then follow a density increment argument, similar to the one described above, to find a rainbow subdivision of $K_k$.

		This bound was subsequently improved to a tighter $O(n (\log n)^{53})$ by Jiang, Letzter, Methuku, and Yepremyan \cite{jiang2021rainbow}. 
		\begin{theorem}[Jiang--Letzter--Methuku--Yepremyan \cite{jiang2021rainbow}] \label{thm:rainbow-subdivision-jlmy}
			Let $k \ge 3$. If $G$ is an $n$-vertex properly edge-coloured graph with no rainbow $K_k$-subdivision, then $e(G) = O(n (\log n)^{53})$.
		\end{theorem}
		Like the improvement of \cite{letzter2023hypergraphs} regarding the Tur\'an number of tight cycles, the key here was to prove that in a given properly-coloured expander there is a rainbow $K_k$-subdivision. To do so, the authors considered almost-regular expanders (meaning that the ratio between the minimum and maximum degree is at most polylogarithmic in $n$), and proved an appropriate existence result for almost-regular expanders. Next, the authors used the fact that random walks mix rapidly in expanders, as well as tools from \cite{janzer2023rainbow} regarding homomorphism counts, to find the desired rainbow $K_k$-subdivision in a sufficiently dense almost-regular expander. The bound in \Cref{thm:rainbow-subdivision-jlmy} was improved further, first by Tomon \cite{tomon2022robust} and then by Wang \cite{wang2022rainbow}; more on this soon.

	\subsection{Cycles with many chords}

		The connection between random walks and almost-regular expanders, described above, was used recently by Dragani\'c, Methuku, Munh\'a-Correia, and Sudakov \cite{draganic2023cycles} to prove the following. 
		\begin{theorem}[Dragani\'c--Methuku--Munh\'a-Correia--Sudakov \cite{draganic2023cycles}] \label{thm:cycle-chords}
			There is a constant $c > 0$ such that every $n$-vertex graph with at least $c \cdot n (\log n)^8$ edges contains a cycle $C$ with at least $|C|$ chords.
		\end{theorem}
		This result is tight up to the polylogarithmic factor. It is plausible that a linear bound in $n$ suffices.
		\begin{question}
			Is there $c > 0$ such that every $n$-vertex graph with at least $cn$ edges contains a cycle $C$ with at least $|C|$ chords?
		\end{question}
		The significantly lower exponent of $\log n$ here, when compared with \Cref{thm:rainbow-subdivision-jlmy}, is at least in part due to the use of a sharper existence result, showing that in a graph $G$ there is an $100$-almost regular (namely, the ration between the maximum and minimum degrees is at most $100$) expander $H$ with $d(H) = \Omega(d(G)/\log n)$.\footnote{the notion of expansion here is similar to the one used in \cite{jiang2021rainbow}.} Using, among other things, the fact that random walks mix rapidly in expanders, the authors show that a random walk $W$ of appropriate length $t$ is, with positive probability: self-avoiding (i.e.\ a path); spans more than $2t$ edges among the middle $t/2$ vertices; and has an edge between the first and last $t/4$ vertices. This clearly yields a cycle of length at most $t$ with at least $t$ chords, as required. 

	\subsection{A sampling trick}

		\Cref{thm:rainbow-subdivision-jlmy} was recently improved by Tomon \cite{tomon2022robust}, who showed that $n (\log n)^{6 + \Omega(1)}$ edges suffice to guarantee the existence of a rainbow $K_k$-subdivision. 
		\begin{theorem}[Tomon \cite{tomon2022robust}] \label{thm:tomon-subdivision}
			Let $k \ge 3$. If $G$ is an $n$-vertex properly edge-coloured graph which has no rainbow $K_k$-subdivisions, then $e(G) \le n \cdot (\log n)^{6+o(1)}$.
		\end{theorem}
		A key component in Tomon's proof is a clever sampling trick, which allowed him to prove the following lemma: 
		\begin{enumerate}[label = \ding{96}]
			\item \label{itm:tomon}
				In a properly edge-coloured expander $H$ (with an appropriate definition of expanders\footnote{Tomon defines an \emph{$\alpha$-maximal} graph to be a graph $G$ satisfying $\frac{d(G)}{|G|^{\alpha}} \ge \frac{d(H)}{|H|^{\alpha}}$ for every subgraph $H \subseteq G$, and shows that $\alpha$-maximal graphs have strong expansion properties. In his paper $\alpha$ is taken to be either $\Theta(1/\log n)$ or a constant, depending on the context.}), if $U$ and $C$ are random sets of vertices and colours, obtained by including each vertex or colour with probability $p$, independently, then for every vertex $v$ at least $\Omega(n)$ vertices in $H$ can be reached from $v$ via a short rainbow path whose colours are in $C$ and interior vertices are in $U$, as long as $d(H)$ is large enough with respect to $p$. 
		\end{enumerate}

		Let us say a few words about the proof of \ref{itm:tomon}. Using the `sprinkling' method, think of $U$ and as the union of smaller disjoint random sets $U_1, \ldots, U_{\ell}$, and similarly think of $C$ as the union of smaller disjoint random $C_1, \ldots, C_{\ell}$. Now define $B_i$ to be the set of vertices reachable from $v$ via a rainbow path of length at most $i$ with interior in $U_1 \cup \ldots \cup U_i$ and colours in $C_1 \cup \ldots \cup C_i$; importantly, vertices in $B_i$ need not be in $U_1 \cup \ldots \cup U_i$. Now it suffices to show that the sets $|B_i|$ grow sufficiently rapidly, until reaching size $\Omega(n)$. To prove this, Tomon distinguishes between the cases where the neighbourhood of $B_i$ has many vertices with few neighbours in $B_i$, and vice versa.

		Lemma \ref{itm:tomon} implies that many pairs of vertices have many rainbow paths joining them, which have pairwise disjoint interiors and pairwise disjoint colour sets. It is then not hard to find a rainbow clique subdivision.
		Variants of this lemma allowed Tomon to improve the bound on the rainbow Tur\'an number of cycles to $\exr(n, \cC) \le n (\log n)^{2 + o(1)}$, and to obtain interesting Tur\'an type results about triangulations of the cylinder and M\"obius strip in $3$-uniform hypergraphs.

	\subsection{Further improvements on rainbow cycles and clique subdivisions}

		Tomon's bound on the number of edges guaranteeing a rainbow $K_t$-subdivision was subsequently improved by Wang \cite{wang2022rainbow}, who showed that $n (\log n)^{2 + \Omega(1)}$ edges suffice.
		\begin{theorem}
			Let $k \ge 3$. If $G$ is an $n$-vertex properly edge-coloured graph with no rainbow $K_k$-subdivisions, then $e(G) \le n \cdot (\log n)^{2 + o(1)}$.
		\end{theorem}

		He did so by optimising Tomon's argument in the context of rainbow clique subdivisions, showing that in a properly edge-coloured expander (Wang used yet another notion which is quite similar to that in \Cref{def:komlos-szemeredi}), for every vertex $v$, if $W$ is a small set of forbidden colours and $C$ is a random colour set obtained by including each colour with probability $1/2$, independently, then, with high probability, more than half the vertices in the graph are reachable from $v$ via a short rainbow path using colours in $C \setminus W$. This readily implies that if $W$ is a small set of forbidden colours, then every two vertices can be joined by a rainbow path. One can then construct a rainbow $K_k$-subdivision greedily.

		Very recently, Tomon's bound on the rainbow Tur\'an number of cycles was improved slightly, to $\exr(n, \cC) = O(n (\log n)^2)$, by Janzer and Sudakov \cite{janzer2022turan}, and, independently, by Kim, Lee, Liu, and Tran \cite{kim2022rainbow}. 
		\begin{theorem}[Janzer--Sudakov \cite{janzer2022turan} and Kim--Lee--Liu--Tran \cite{kim2022rainbow}]
			Let $G$ be an $n$-vertex properly edge-coloured graph with no rainbow cycles. Then $e(G) = O(n \cdot (\log n)^2)$.
		\end{theorem}

		Neither paper used expanders. Instead, they used inequalities regarding homomorphism counts, reminiscent of Janzer's methods from \cite{janzer2023rainbow}, with \cite{janzer2022turan} considering a weighted count, and \cite{kim2022rainbow} using an unweighted count along with regularisation.

	\subsection{An almost tight result regarding rainbow cycles}

		Even more recently\footnote{and after this survey was submitted for publication}, Alon, Buci\'c, Sauermann, Zakharov, and Zamir \cite{alon2023essentially} proved the following almost tight bound on the number of edges in a properly coloured graph with no rainbow cycles.
		\begin{theorem}[Alon--Buci\'c--Sauermann--Zakharov--Zamir \cite{alon2023essentially}] \label{thm:rainbow-cycle-almost-tight}
			Let $G$ be a properly edge-coloured $n$-vertex graph with no rainbow cycles. Then $e(G) \le O(n \log n \log \log n)$.
		\end{theorem}

		In their proof of \Cref{thm:rainbow-cycle-almost-tight} the authors use the following notion of expanders.

		\begin{definition}[Alon--Buci\'c--Sauermann--Zakharov--Zamir \cite{alon2023essentially}] \label{def:ab-expander}
			An $n$-vertex graph $G$ is an \emph{\abExpander} if, for every $\eps$ with $0 \le \eps \le 1$, every subset $U \subseteq V(G)$ with $1 \le |U| \le n^{1-\eps}$, and every subset $F \subseteq E(G)$ with $|F| \le (\eps/3) d(G) |U|$, we have $|N_{G - F}(U)| \ge (\eps/3)|U|$.
		\end{definition}

		Observe that in an $n$-vertex \abExpander, vertex sets $U$ of size at most $n^{0.99}$, say, expand linearly (namely $|N(U)| \ge \frac{1}{100}|U|$, because we can take $\eps = \frac{1}{100}$), and set of size at most $n/2$ expand at a rate of at least $\frac{1}{\log n}$ (by taking $\eps = \frac{1}{\log n}$).

		These expanders are relevant in the setting of \Cref{thm:rainbow-cycle-almost-tight} due to the following lemma.
		\begin{lemma}[Alon--Buci\'c--Sauermann--Zakharov--Zamir \cite{alon2023essentially}] \label{lem:ab-expander-existence}
			Every graph $G$ with at least one edge contains an \abExpander $H$ with 
			\begin{equation*}
				d(H) \ge \frac{1}{3} \cdot \frac{\log |H|}{\log |G|} \cdot d(G).
			\end{equation*}
		\end{lemma}

		In particular, this lemma implies that if $d(G) \ge 3c \cdot \log |G| \log \log |G|$ then there is an \abExpander $H \subseteq G$ with $d(H) \ge c \cdot \log |H| \log \log |H|$.
		Thus, it suffices to show that, if $H$ is a properly coloured $n$-vertex \abExpander with $d(H) \ge c \cdot \log n \log \log n$, where $c$ is a large constant, then $H$ has a rainbow cycle.
		Similarly to Tomon \cite{tomon2022robust} (see \ref{itm:tomon}), the authors prove the following.
		\begin{enumerate}[label = \ding{93}]
			\item \label{itm:abszz}
				Let $H$ be an $n$-vertex properly coloured \abExpander with $d(H) \ge c \cdot \log n \log \log n$, where $c$ is a large constant, and let $v$ be a vertex in $H$.
				Let $C$ be a random set of colours, obtained by including each colour in $H$ with probability $1/2$, independently.
				Then, with probability larger than $1/2$, at least $\frac{n+1}{2}$ vertices in $H$ can be reached from $v$ via a rainbow path whose colours are in $C$.
		\end{enumerate}

		Notice that \ref{itm:abszz} readily implies \Cref{thm:rainbow-cycle-almost-tight}. Indeed, as pointed out, it suffices to show that every expander $H$ as in \ref{itm:abszz} contains a rainbow cycle. Let $\{C_1, C_2\}$ be a random partition of the colour set of $H$. Then, by \ref{itm:abszz}, with positive probability the sets $U_1, U_2$, of vertices that are reachable from $v$ through a rainbow path with colours in $C_1, C_2$, respectively, satisfy $|U_1|, |U_2| > n/2$. Fix such a partition $\{C_1, C_2\}$, and let $u \in U_1 \cap U_2$. Then there is a rainbow closed walk through $v$ and $u$, which in turn contains a rainbow cycle, as required.

		A key novelty in the proof of \ref{itm:abszz} is another sampling trick.
		Define $C_0$ as above, namely to include each colour in $H$ with probability $1/2$, independently.
		For appropriate $t$ and $p$, define $C_1, \ldots, C_t$ so that $C_i$ is obtained from $C_{i-1}$ by removing each element of $C_{i-1}$ with probability $p$, independently. Now define $U_i$ to be the set of vertices reachable from $v$ via a rainbow path with colours in $C_i$.
		Roughly speaking, the authors of \cite{alon2023essentially} show that, in expectation, 
		\begin{equation} \label{eqn:Ui}
			|U_{i-1}| \ge (1 + \Omega(\eps))|U_{i+1}|,
		\end{equation}
		where $\eps$ satisfies $|U_i| = n^{1-\eps}$.
		Assuming the $U_i$'s indeed follow this behaviour (not just in expectation), then the parameters are such that $|U_0| > n/2$, as required.
		Turning this into an actual proof involves various probabilistic ideas that we do not mention here. Instead, we sketch how to prove that \eqref{eqn:Ui} holds in expectation. 

		Fix $U_i$, a corresponding $\eps$ (so that $|U_i| = n^{1-\eps}$), and $C_i$. Then every colour in $C_i$ is not in $C_{i+1}$ with probability $p$, and one can check that every colour not in $C_i$ is in $C_{i-1}$ with probability $\Omega(p)$.
		The authors first prove a dichotomy, somewhat in spirit of \ref{itm:bm-connecting} below, according to which there are many edges $E$ leaving $U_i$ such that either $E$ is coloured in $C_i$ and no vertex in $U_i$ is incident to many edges of $E$, or no edges in $E$ are coloured in $C_i$ and no vertex outside of $U_i$ is incident to many edges in $E$.

		Suppose the latter holds. Notice that every vertex outside $U_i$ that is incident to an edge in $E$ whose colour is in $C_{i-1}$ is in $U_{i-1}$. Hence, using that $E$ is large, no vertex outside of $U_i$ is incident with many edges in $E$, and that the colour of every edge in $E$ is in $C_{i-1}$ with probability $\Omega(p)$, we get that, in expectation, $|U_{i-1}| \ge (1 + \Omega(\eps))|U_i|$.
		Now consider the former case. Observe that if $u \in U_i$ is incident to an edge $uw$ in $E$ whose colour is not in $C_{i+1}$ then $u \notin U_{i+1}$. Indeed, otherwise $w$ would be in $U_i$ (because there would be a rainbow path from $v$ to $u$ coloured in $C_{i+1}$, which can be extended by $uw$ to a rainbow path from $v$ to $w$ coloured in $C_{i}$), a contradiction. Similar arguments to the previous case now show that, in expectation, $|U_i| \ge (1 + \Omega(\eps))|U_{i+1}|$.
		Either way, \eqref{itm:abszz} holds.

	\subsection{Cycles with all diagonals}

		A \emph{diagonal} in a cycle $C$ is a chord joining two vertices at distance $\floor{|C|/2}$ in $C$.
		Erd\H{o}s \cite{erdos1975some} asked the following question regarding cycles with all diagonals: what is the maximum number of edges in a graph on $n$ vertices that has no cycle containing all diagonals? 
		\begin{figure}[h]
			\centering
			\includegraphics[scale = .7]{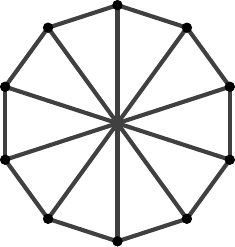}
			\caption{A cycle with all diagonals}
			\label{fig:cycle-diagonals}
		\end{figure}
		Notice that a cycle with all diagonals contains every cycle of length $4$. Hence, since there are graphs on $n$ vertices with $\Omega(n^{3/2})$ edges and no $4$-cycles, the answer is at least $\Omega(n^{3/2})$.
		Erd\H{o}s observed also that $K_{3,3}$ is a $6$-cycle with all diagonals, and so the answer is at most the Tur\'an number $\ex(n, K_{3,3})$ of $K_{3,3}$, which is $O(n^{5/3})$.
		Very recently, Brada\'c, Methuku, and Sudakov \cite{bradavc2023extremal} gave a tight answer to the above question (up to a constant factor).
		\begin{theorem}[Brada\'c--Methuku--Sudakov \cite{bradavc2023extremal}] \label{thm:cycle-diagonals}
			Every $n$-vertex graph with no cycles containing all diagonals has $O(n^{3/2})$ edges.
		\end{theorem}

		While the proof of \Cref{thm:cycle-diagonals} uses expanders, it does not fall under the scope of this survey, as the expansion rate in their expanders is linear. Nevertheless, their proof uses the notion of expansion introduced by Tomon \cite{tomon2022robust} in his proof of \Cref{thm:tomon-subdivision}, and ideas of the author \cite{letzter2023hypergraphs} used in her proof of \Cref{thm:hypergraph}, so we briefly sketch it. 

\begin{figure}[h]
			\centering
			\includegraphics[scale = .6]{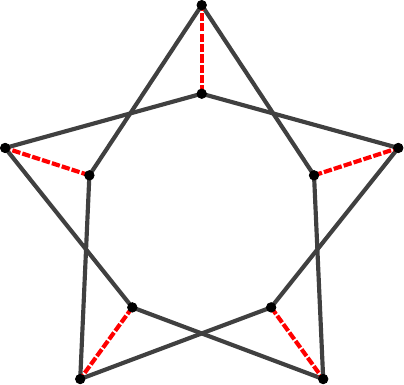}
			\caption{A $5$-cycle in $\Gamma$ whose vertices are pairwise vertex-disjoint edges in $H$ (denoted by red striped edges). The vertices of the $5$-cycle are the diagonals in a $10$-cycle in $H$.}
			\label{fig:cycle-diagonals-var}
		\end{figure}

		An interesting new idea here is a way to find an almost spanning expander.
		They start with a bipartite $n$-vertex graph with at least $cn^{3/2}$ edges, for some (large) constant $c$. Within this graph, they find an expander (using Tomon's setup and an additional cleanup step) $H$, which has $m$ vertices, at least $c' m^{3/2}$ edges for an appropriate constant $c'$, and maximum degree $O(m^{1/2})$.
		They then define an auxiliary graph $\Gamma$, whose vertices are the edges of $H$, and where two edges $xy$ and $uv$ are joined if $(xyuv)$ is a $4$-cycle in $H$.
		The main novelty in their proof is a way of finding an expander $\Gamma'$ which is an almost spanning subgraph of $\Gamma$. This is achieved by showing that every set of vertices in $\Gamma$, which is not too large or too small, expands.
		They then use the expansion properties of $\Gamma'$, and ideas from \cite{letzter2023hypergraphs}, to find an odd cycle $(e_1 \ldots e_{\ell})$ in $\Gamma$, where the $e_i$'s, viewed as edges of $H$, are pairwise vertex-disjoint. This readily implies the existence of a $2\ell$-cycle in $H$ that contains all diagonals (using that $H$ is bipartite; see \Cref{fig:cycle-diagonals-var}).

\section{Decompositions and separation} \label{sec:decompositions-separation}

	All the results we have mentioned so far find certain substructures in graphs $G$ with some properties (typically, with large average degree), which might be much smaller or sparse than $G$. In particular, these results do not aspire to say anything global about $G$. Often, they find an expander $H$ within $G$, which might be much smaller than $G$, and focus just on $H$.
	In this section we will see how sublinear expanders can be used to prove more global results.

	\subsection{Decomposing graphs into cycles and edges}

		A \emph{decomposition} of a graph $G$ is a collection of subgraphs of $G$ such that each edge in $G$ appears in exactly one of these subgraphs.
		In the 1960's a lot of questions regarding decompositions of a given graph into simpler subgraphs were considered (see, e.g., \cite{lovasz1968covering,harary1970covering,vizing1964estimate}). One example of a conjecture along these lines is the following conjecture due to Erd\H{o}s and Gallai. 
		\begin{conjecture}[Erd\H{o}s--Gallai (see \cite{erdos1966representation})] \label{conj:erdos-gallai}
			There is a constant $c > 0$ such that every $n$-vertex graph can be decomposed into at most $cn$ cycles and edges.
		\end{conjecture}
		Notice that allowing edges in the decomposition is necessary, due to the potential existence of odd degree vertices. Moreover, at least $n-1$ edges and cycles are necessary, as can be seen by considering a tree.
		A greedy algorithm that picks a longest cycle and removes it from the graph, one by one, yields a decomposition into $O(n \log n)$ cycles and edges. The first improvement on this was achieved in 2014 by Conlon, Fox, and Sudakov \cite{conlon2014cycle}, who proved that $O(n \log \log n)$ cycles and edges suffice. This was recently improved upon significantly by Buci\'c and Montgomery \cite{bucic2022towards}.

		\begin{theorem}[Buci\'c--Montgomery \cite{bucic2022towards}] \label{thm:erdos-gallai}
			Every $n$-vertex graph can be decomposed into $O(n \logstar n)$ cycles and edges.
		\end{theorem}

		Here $\logstar n$ is the \emph{iterated logarithm}, namely the minimum number of times the (say) base $2$ logarithm has to be applied, sequentially, to $n$, to reach a number which is at most $1$.

		To prove \Cref{thm:erdos-gallai}, the authors iterate the following statement: 
		\begin{enumerate}[label = \ding{118}]
			\item \label{itm:bm-statement}
				If $G$ is an $n$-vertex graph with average degree $d$, then $G$ can be decomposed into $O(n)$ cycles and $O\big(n (\log d)^{O(1)}\big)$ edges.
		\end{enumerate}
		Notice that \ref{itm:bm-statement} implies \Cref{thm:erdos-gallai}: iterating this statement $c \logstar n$ times, for large enough $c$, brings the average degree down from at most $n$ to at most $O(1)$, using $O(n)$ edges at each step, and then the remaining edges can be covered by single edges, thereby using $O(n \logstar n)$ edges in total.

		The high level sketch of the proof of \ref{itm:bm-statement} is as follows. 
		\begin{itemize}
			\item
				First, remove, one by one and as long as possible, cycles of length at least $d$.
			\item
				Next, decompose the remainder into (sublinear) expanders (on at most $O\big(d (\log d)^{O(1)}\big)$ vertices), which are almost vertex-disjoint.
			\item
				Finally, decompose each expander $H$ into $O(|H|)$ cycles and $O\big(|H|(\log d)^{O(1)}\big)$ edges.
		\end{itemize}
		This outline is inspired by the proof of Conlon, Fox, and Sudakov \cite{conlon2014cycle}, who used much denser expanders, to show that $O(n)$ cycles suffice to reduce the average degree from $d$ to $O(d^{1-\Omega(1)})$. As such, the details are rather different, and quite a lot of new ideas were needed here. 
		
		The definition of expanders used in \cite{bucic2022towards} is similar to the robust expanders defined in \Cref{def:komlos-szemeredi-robust}.
		\begin{definition}[Buci\'c--Montgomery \cite{bucic2022towards}]
			An $n$-vertex graph $G$ is an \emph{$(\eps,s)$-\bmExpander} if for every vertex set $U \subseteq V(G)$ and subgraph $F \subseteq G$, satisfying $|U| \le 2n/3$ and $e(F) \le s|U|$, the following holds: $|N_{G \setminus F}(U)| \ge \frac{\eps|U|}{(\log n)^2}$.
		\end{definition}

		To realise the second step above, the authors prove a decomposition result. 
		\begin{lemma}[Lemma 14 in \cite{bucic2022towards}] \label{lem:decomposition}
			Let $G$ be an $n$-vertex graph. Then it can be decomposed into $(\eps,s)$-\bmExpanders $G_1, \ldots, G_r$, with $\sum|G_i| = O(n)$, and $O(sn\log n)$ edges.
		\end{lemma}
		To complete the second step, the authors apply the lemma with $s = 0$, and prove that every $n$-vertex $(\eps,0)$-expander contains a cycle of length $\Omega(n (\log n)^{-4})$, using a DFS algorithm. This shows that the expanders in the decomposition are rather small (they have at most $O\big( d (\log d)^{O(1)}\big)$ vertices), because the graph that remains after the first step has no cycles of length at least $d$.

		The main challenge is in the proof of the final step. Applying the \Cref{lem:decomposition} again to each expander $H$, with $s = (\log n)^c$ for some constant $c$, it suffices to show that every $n$-vertex $(\eps, s)$-expander can be decomposed into $O(n)$ cycles and at most $n (\log n)^{O(1)}$ edges.
		
		A key result in this direction is the following (see Lemma 19 in \cite{bucic2022towards}). 
		\begin{enumerate}[label = \ding{71}]
			\item \label{itm:bm-connecting}
				Given an $n$-vertex $(\eps,s)$-expander $H$, if $V$ is a random set of vertices that includes each vertex with probability $1/3$, independently, then, with high probability, for every vertex set $U$, and small subgraph $F$ (of size at most $|U|/(\log n)^{27}$), more that $|V|/2$ vertices in $V$ can be reached from $U$ through a relatively short path with interior in $V$.
		\end{enumerate}

		The proof of this uses the sampling trick of Tomon \cite{tomon2022robust} described before (see \Cref{thm:tomon-subdivision} and \ref{itm:tomon}), as well as an expansion dichotomy for vertex sets $U$ in $(\eps,s)$-expanders (see Propositions 12 and 13 in \cite{bucic2022towards}): 
		\begin{enumerate}[label = \ding{78}]
			\item
				If $|U| \le 2n/3$ and $e(F) \le s|U|/2$, then either the neighbourhood $N_{G \setminus F}(U)$ is very large, or many vertices in the neighbourhood $N_{G \setminus F}(U)$ have many neighbours in $U$.
		\end{enumerate}
		The authors use this dichotomy to apply the union bound effectively.
		Using the key result \ref{itm:bm-connecting}, and a result about matchings in hypergraphs due to Aharoni and Haxell \cite{aharoni2000hall}, they prove two results about joining pairs of vertices through a random vertex set, which they then use to complete the third step.

	\subsection{Separating the edges of a graph by paths}

		A \emph{separating path system} for a graph $G$ is a collection $\cP$ of paths such that for every two edges $e$ and $f$ there is a path $P \in \cP$ that contains $e$ but not $f$.
		This notion was introduced by Katona (2013). Writing $\sep(G)$ for the size of a smallest separating path system for $G$, and $\sep(n)$ for the maximum of $\sep(G)$ over all $n$-vertex graphs, Katona asked to determine $\sep(n)$.
		We claim that $\sep(n) = O(n \log n)$. Indeed, given an $n$-vertex graph $G$, consider a collection $\cG$ of $O(\log n)$ subgraphs of $G$ that separates the edges of $G$, meaning that for every two edges $e$ and $f$ in $G$ there is a subgraph in $\cG$ that contains $e$ but not $f$ (that such a collection exists follows by noticing that every set of size $m$ can be separated using $O(\log m)$ sets, an easy exercise).
		Using a result of Lov\'asz \cite{lovasz1968covering}, asserting that the edges of every $n$-vertex graph can be decomposed into at most $n$ paths, for each $H \in \cG$ there is a collection $\cP_H$ of at most $n$ paths decomposing the edges of $H$. 
		It is easy to verify that $\bigcup_{H \in \cG}\cP_H$ is a separating path system for $G$, which has size $O(n \log n)$.

		The first improvement on this initial bound was achieved by the author \cite{letzter2022separating}.
		\begin{theorem}[Letzter \cite{letzter2022separating}] \label{thm:sep}
			Every $n$-vertex graph has a separating path system of size $O(n \logstar n)$.
		\end{theorem}

		The proof draws on many ideas from Buci\'c and Montgomery \cite{bucic2022towards}.
		One difference is that here, unlike in \cite{bucic2022towards}, one needs to work with relatively sparse graphs, making the decomposition lemma (\Cref{lem:decomposition}) too weak at times. To overcome this, we introduce a variant of \bmExpanders.
		\begin{definition}[Letzter \cite{letzter2022separating}]
			An $n$-vertex graph $G$ is an $(\eps, s, t)$-\lExpander if for every vertex set $U \subseteq V(G)$ and subgraph $F \subseteq G$, such that $1 \le |U| \le 2n/3$ and $|F| \le s \cdot \min\{|U|, t\}$, we have $|N_{G \setminus F}(U)| \ge \frac{\eps|U|}{(\log |U| + 1)^2}$.
		\end{definition}
		We also prove an analogue of \Cref{lem:decomposition}.
		\begin{lemma} \label{lem:decomposition-sep}
			If $G$ is an $n$-vertex graph, then it can be decomposed into $(\eps, s, t)$-\lExpanders $G_1, \ldots, G_r$ with $\sum |G_i| = O(n)$ and $O(sn (\log t)^2)$ edges.
		\end{lemma}

		The role of the parameter $t$ is to allow for the expanders in the decomposition to have good expansion properties for small sets $U$, while not causing too many edges in $G$ to remain uncovered.

		To prove \Cref{thm:sep}, a basic idea is that it is easy to find, given an $n$-vertex graph $G$, collections of paths $\cP$ and matchings $\cM$, of size $O(n)$, that separate the edges of $G$. As such, using the decomposition lemma, it suffices to be able to extend a given matching $M$ to a path that avoids a small subgraph (corresponding to paths in $\cP$ that intersect $M$). This is achieved for dense expanders using variants of tools from \cite{bucic2022towards}, and for sparse graphs using methods developed for \ksExpanders, dealing separately with edges touching large degree vertices.

		Shortly after \cite{letzter2022separating} appeared, Bonamy, Botler, Dross, Naia, and Skokan \cite{bonamy2023separating} proved that $\sep(n) \le 19n$, using a simple inductive argument. This is tight up to the factor 19, confirming a conjecture from \cite{balogh2016path} and \cite{falgas2013separating}. They raised the following question.

		\begin{question}[Bonamy--Botler--Dross--Naia--Skokan \cite{bonamy2023separating}] \label{qn:rainbow-path-cover}
			Is there a constant $c > 0$ such that the edges of every properly edge-coloured $n$-vertex graph can be covered by $O(n)$ rainbow paths?
		\end{question}

		This can be shown to imply the existence of a separating path system of size $O(n)$, via Lov\'asz's result about decomposing a graph into paths. 

\section{Counting Hamiltonian sets} \label{sec:hamiltonian}

	An interesting different direction that has been addressed using sublinear expander involves the number of \emph{Hamiltonian sets}. In a graph $G$, a set of vertices $U$ is called \emph{Hamiltonian} if $G[U]$ has a Hamiltonian cycle. Denote by $h(G)$ the number of Hamiltonian sets in $G$. 

	\subsection{Maximising $h(G)$ among graphs with given average degree}
		Koml\'os (see \cite{tuza1990exponentially,tuza2001unsolved,tuza2013problems}) conjectured that among graphs $G$ with minimum degree at least $d$, the complete graph $K_{d+1}$ minimises $h(G)$; namely: if $\delta(G) \ge d$ then $h(G) \ge h(K_{d+1})$.
		In 2017, Kim, Liu, Sharifzadeh, and Staden \cite{kim2017proof} proved this conjecture (for large $d$).
		\begin{theorem}[Kim--Liu--Sharifzadeh--Staden \cite{kim2017proof}] \label{thm:komlos-conj}
			Let $d$ be large, and suppose that $G$ is a graph with $d(G) \ge d$. Then $G$ has at least $h(K_{d+1}) = 2^{d+1} - \binom{d+1}{2} - (d+1) - 1$ Hamiltonian sets.
		\end{theorem}
		In fact, they prove that if $d(G) \ge d$ and $G$ is not isomorphic to $K_{d+1}$ or to the union of two $K_{d+1}$'s that share a single vertex, then $h(G) \ge (2 + o(1))2^{d+1}$ (notice that $h(K_{d+1}) \approx 2^{d+1}$ and $h(H) \approx \frac{3}{2} \cdot 2^{d+1}$ for $H$ the union of two $K_{d+1}$'s sharing a single vertex).
		They also sketch how their methods can be used to address a bipartite version of the same problem.

		As a first step, the authors find an $n$-vertex $(\eps_1, \eps_2d)$-expander $H$ with $\delta(H) = \Omega(d)$ and $\Delta(H) = O(d)$, with $d = o(n)$.
		To do so, suppose there is no such expander. By analysis of ``blocks'' (maximal $2$-connected subgraphs), assuming that $G$ is a minimal counterexample to \Cref{thm:komlos-conj}, and by a proof of \Cref{thm:komlos-conj} for dense graphs (using the regularity lemma), they may essentially assume that $G$ is $2$-connected and of order $\omega(d)$. Now, after removing the few large degree vertices, either there are two disjoint dense subexpanders $H_1,H_2$ with average degree a at least a bit below $d$, or there is a subexpander $H$ with the required properties. Either way, this leads to a contradiction: either to the $G$ being a counterexample, or to there being no expanders as above.

		Now, given an almost-regular expander $H$ as above, the authors again consider two cases: dense ($d \ge (\log n)^{c}$) and sparse (otherwise).
		In the sparse case, they find a set $Z$ of $200d$ vertices that are far apart, and show that for every $U \subseteq Z$ of size $100d$ there is a cycle $C$ with $V(C) \cap Z = U$, yielding $\binom{200d}{100d}$ Hamiltonian sets. This can be done in a routine way, by joining pairs of vertices in $U$ by short paths avoiding the vicinity of other vertices in $U$.
		The dense case is a little more involved. They first find $200d$ many $(h_0,h_1,h_2,h_3)$-webs (see \Cref{def:web}; here $h_0,h_1,h_3 = (\log n)^{\Theta(1)}$ and $h_2 = \Theta(d)$), whose sets of non-leaf vertices are pairwise disjoint. Denoting the set of roots by $Z$, for every $U \subseteq Z$ of size $100d$, they find a cycle $C$ that avoids $Z \setminus U$ and contains at least $98d$ vertices of $U$ (this yields at least the following number of different Hamiltonian sets $\binom{200d}{100d} / \binom{102d}{2d} \ge \Theta(d^{-1/2} 2^{200d} / 2^{102d}) \ge 2^{50d}$, with room to spare in the last inequality). To do so, they join vertices in $U$ one by one through the webs, avoiding vertices with overused leaf sets.

	\subsection{Maximising $h(G)$ among graphs with given average degree and number of vertices}

		Cambie, Gao, and Liu \cite{cambie2023many} considered a similar question: among $n$-vertex graphs $G$ with $\delta(G) \ge d$, how small can $h(G)$ be? To describe their result, we need the notion of `crux', given in \Cref{def:crux}.
		\begin{theorem}[Cambie--Gao--Liu \cite{cambie2023many}] \label{thm:hamiltonian-sets}
			Let $G$ be a graph, $d = d(G)$, and $t = c_{1/5}(G)$.
			Then 
			\begin{equation*}
				h(G) \ge n \cdot 2^{\Omega(t\, (\log t)^{-16})}.
			\end{equation*}
		\end{theorem}
		At least in some regimes, this is tight up to the polylogarithmic factor, as can be seen by considering a disjoint union of $\floor{n/(d+1)}$ copies of $K_{d+1}$.
		
\begin{figure}[ht]
			\centering
			\includegraphics[scale = .7]{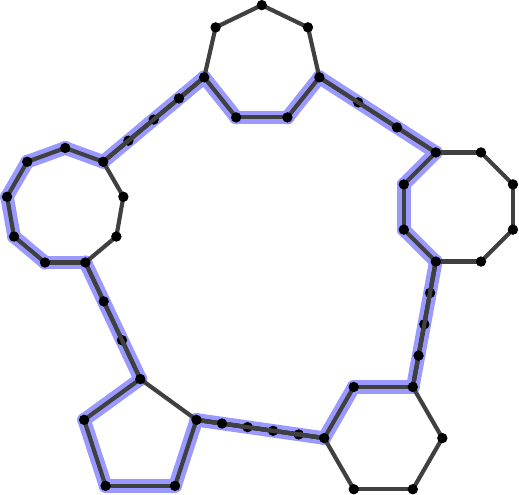}
			\caption{A `chain' of cycles. Each of the vertices in the `connecting paths' is in $32$ cycles; one of them is marked in blue.}
			\label{fig:chain}
		\end{figure}

		The main technical result towards the proof of the above is the following. 
		\begin{enumerate}[label = \ding{68}]
			\item \label{itm:hamiltonian}
				If $|G| = n$, $d(G) \ge d$, and $t = c_{\alpha}(G)$, then there is a vertex in at least $2^{\Omega(t\,(\log t)^{-c})}$ different Hamiltonian sets.
		\end{enumerate}
		To prove \ref{itm:hamiltonian}, they consider an $(\eps_1, \eps_2d)$-\ksExpander $H$ on $k$ vertices, find many disjoint cycles of length $(\log k)^{\Theta(1)}$, and join $k (\log k)^{-\Theta(1)}$ into a `chain' by vertex-disjoint paths (see \Cref{fig:chain}). Then, every vertex in one of the connecting paths is in $2^{k(\log k)^{-\Theta(1)}}$ different Hamiltonian sets, because in forming a cycle through the chain, for each cycle there are two possible choices of a path through it to the next connecting path. By the definition of crux, we have $k \ge t$, and hence there is a vertex in at least the required number of different Hamiltonian sets.

		Now, to get the main result, consider a graph $G$ which is a minimal counterexample to the main result. First, remove, one by one, vertices which are in many different Hamiltonian sets, and denote the sets of removed vertices by $S$. Notice that $G \setminus S$ has small average degree, because otherwise one could apply \ref{itm:hamiltonian} to a subexpander of $G \setminus S$ with large enough average degree, contradicting the choice of $S$.
		Next, repeat a similar procedure in $G[S, V(G) \setminus S]$, and move, one by one, a vertex from $V(G) \setminus S$ to $S$ if if is in many different Hamiltonian sets in $G[S, V(G) \setminus S]$. If the resulting $S$ is large, then $h(G)$ is large enough for a contradiction. Otherwise, it is small but has large average degree. This gives an upper bound on the crux of $G$, which together with the minimality assumption on $G$ yields a contradiction.

		Recall that \Cref{thm:hamiltonian-sets} is close to tight, but unlikely to be tight. As such, it might be interesting to try to improve the bound.
		\begin{problem}
			Prove a tighter lower bound on $h(G)$ for an $n$-vertex graph $G$ with average degree at least $d$. In particular, is it true that $h(G) \ge n \cdot 2^{\Omega(t)}$, where $t = c_{1/5}(G)$?
		\end{problem}

\section{Other results} \label{sec:other} 

	In this section we briefly mention a few other results, that used sublinear expanders in their proof.
	The first two are extremal results, showing that a large average degree implies the existence of a certain structure, the third is a result in Ramsey theory.

	\subsection{Nested cycles}
\begin{figure}[ht]
			\centering
			\begin{subfigure}[b]{.3\textwidth}
				\centering
				\includegraphics[scale = .8]{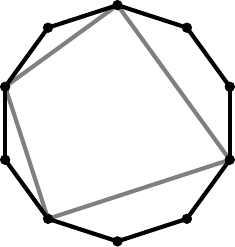}
			\end{subfigure}
			\hspace{1cm}
			\begin{subfigure}[b]{.3\textwidth}
				\centering
				\includegraphics[scale = .8]{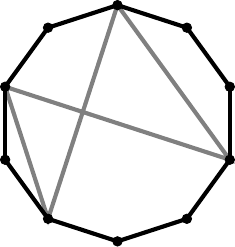}
			\end{subfigure}
			\caption{Two pairs of nested cycles (without and with crossings)}
			\label{fig:nested}
		\end{figure}

		In 1975 Erd\H{o}s \cite{erdos1975some} made two conjectures about the existence of nested cycles in graph with large average degree. The first conjecture asserts that there exists $d_0$ such that every graph $G$ with $d(G) \ge d_0$ contains two edge-disjoint cycles $C_1, C_2$ such that $V(C_2) \subseteq V(C_1)$. The second is a strengthening of the former, additionally requiring that the cyclic ordering of $C_2$ respects that of $C_1$. More precisely, if $C_1 = (v_1 \ldots v_{\ell})$ then there exist $i_1, \ldots, i_k$, such that $1 \le i_1 < \ldots < i_k \le \ell$ and $C_2 = (v_{i_1} \ldots v_{i_k})$. We refer to a sequence of nested cycles $C_1, \ldots, C_k$ (with $V(C_1) \supseteq \ldots \supseteq V(C_k)$), where the cyclic order of $C_{i+1}$ respects that of $C_i$, as a sequence of \emph{nested cycles without crossings} (see \Cref{fig:nested}). 

		The former conjecture was proved by Bollob\'as \cite{bollobas1978nested} in 1978, and extended by Chen, Erd\H{o}s, and Staton \cite{chen1996proof} in 1996 to longer sequences of nested cycles\footnote{Bollob\'as \cite{bollobas1978nested} proved that an average degree of 14 guarantees the existence of two edge-disjoint nested cycles, and this was improved by Chen, Erd\H{o}s, and Staton \cite{chen1996proof} to 10. For general $k$ they showed that an average degree of $6^k$ guarantees the existence of a sequence of $k$ nested edge-disjoint cycles.}.
		Gil Fern\'andez, Kim, Kim, and Liu \cite{gilfernandez2022nested} recently proved the latter conjecture.
		\begin{theorem}[Gil Fern\'andez--Kim--Kim--Liu \cite{gilfernandez2022nested}] \label{thm:nested}
			There exists $d_0$ such that every graph $G$ with $d(G) \ge d_0$ has two nested cycles without crossings.
		\end{theorem}
		The authors do not give an explicit value of $d_0$ for which the statement holds. Roughly speaking, their proof yields a value of $d_0$ which is large, but not regularity lemma large. I am not aware of any non-trivial lower bounds on $d_0$ for which the statement in \Cref{thm:nested} holds.

		Very briefly, the idea here is to take a shortest cycle $C = (v_1 \ldots v_{\ell})$ in an expander $H$, which is going to be the shorter cycle $C_2$, and then it suffices to find paths $P_1, \ldots, P_{\ell}$, whose interiors are pairwise vertex-disjoint and vertex-disjoint of $C$, such that $P_i$ joins $v_i$ with $v_{i+1}$ (where $v_{\ell+1} := v_1$). To find these paths, they first join each vertex $v_i$ to either two vertices of large degree, or to two large sets with small diameter (they call this structure a `kraken'). Next, they join up these vertices of large degree, or the large sets with small diameter, to obtain the required connections between the vertices. 

		\Cref{thm:nested} raises the following interesting question.

		\begin{question} \label{qn:nested}
			Does there exist $d_0(k)$ such that, for every $k \ge 3$, every graph $G$ with $d(G) \ge d_0(k)$ has a sequence $C_1, \ldots, C_k$ of nested cycles with no crossings?
		\end{question}

	\subsection{Pillars}

\begin{figure}[ht]
			\centering
			\includegraphics[scale = .6]{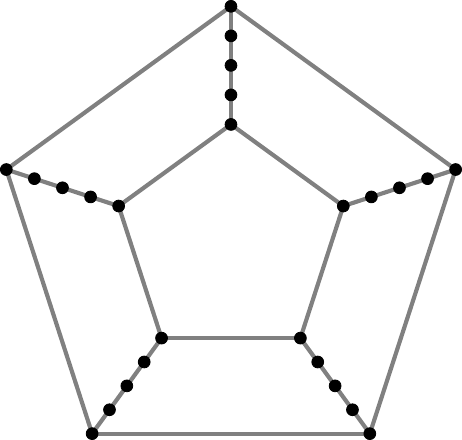}
			\caption{A pillar}
			\label{fig:pillar}
		\end{figure}

		In 1989 Thomassen \cite{thomassen1989configurations} conjectured that large average degree guarantees the existence of a `pillar', defined to consist of two vertex-disjoint cycles of the same length, denoted $C_1 = (v_1 \ldots v_{\ell})$ and $C_2 = (u_1 \ldots u_{\ell})$, and vertex-disjoint paths $P_1, \ldots, P_{\ell}$ of the same length such that $P_i$ joins $v_i$ with $u_i$ (see \Cref{fig:pillar}).

		His conjecture was confirmed by Gil Fern\'andez and Liu \cite{gilfernandez2023how}.
		\begin{theorem}[Gil Fern\'andez and Liu \cite{gilfernandez2023how}] \label{thm:pillar}
			There exists $d_0$ such that every graph $G$ with $d(G) \ge d_0$ contains a pillar.
		\end{theorem}
		In fact, the authors proved a more general result, about the existence of `$K_k$-pillars' (involving $k$ cycles, any two of which are joined as above).

		To prove this, the authors work with an expander $H$ which is $Q_3$-free (because the hypercube $Q_3$ contains a pillar). They then find many disjoint large `krakens' (the structures used in \cite{gilfernandez2022nested}). Two of these would have the same cycle length, and then the kraken's `legs' can be joined together by disjoint paths of the same length, using the approach of Liu and Montgomery \cite{liu2023solution}.

	\subsection{Ramsey goodness of cycles}

		A final, quite different result we mention here involves Ramsey numbers.
		The \emph{Ramsey number} of $H_1$ and $H_2$, denoted $r(H_1, H_2)$, is the minimum $n$ such that every red-blue colouring of $K_n$ contains a red copy of $H_1$ or a blue copy of $H_2$.
		Here we consider the case where $H_1$ is a cycle.
		It is well-known (and easy to see) that 
		\begin{equation} \label{eqn:ramsey-good}
			r(C_n, H) \ge (n-1)(\chi(H) - 1) + \sigma(H),
		\end{equation}
		where $\chi(H)$ is the chromatic number of $H$, and $\sigma(H)$ is the minimum possible size of a colour class in a proper $\chi(H)$-colouring of $H$. 
		Burr \cite{burr1981ramsey} proved that, if $n$ is sufficiently large in terms of $H$, then the bound in \eqref{eqn:ramsey-good} is tight (in which case it is said that $C_n$ is \emph{$H$-good}). Allen, Brightwell, and Skokan \cite{allen2013ramsey} conjectured that the bound is tight already when $n \ge |H| \cdot \chi(H)$.
		Haslegrave, Hyde, Kim, and Liu \cite{haslegrave2023ramsey} proved an even stronger statement, as follows. 
		\begin{theorem}[Haslegrave--Hyde--Kim--Liu \cite{haslegrave2023ramsey}] \label{thm:ramsey-good}
			There exists a constant $c > 0$ such that, if $n \ge c|H|(\log \chi(H))^4$, then $r(C_n, H) = (n-1)(\chi(H)-1) + \sigma(H)$.
		\end{theorem}
		Their approach uses ideas from Liu and Montgomery \cite{liu2023solution}, specifically about the construction of adjusters (see \Cref{def:adjusters}) in sublinear expanders, and is also inspired by a result of Pokrovskiy and Sudakov \cite{pokrovskiy2020ramsey} who proved that \eqref{eqn:ramsey-good} is tight when $n \ge 10^{60}|H|$ and $\sigma(H) \ge \chi(H)^{22}$.

		The lower bound on $n$ in \Cref{thm:ramsey-good} is close to optimal; the authors provide an example showing that $n$ has to be taken to be at least $|H|(1 + o(1))$.
		It would be nice to get rid of the polylogarithmic factor in $\chi(H)$.
		\begin{question}[Haslegrave--Hyde--Kim--Liu \cite{haslegrave2023ramsey}]
			Is there a constant $c > 0$ such that, if $n \ge c|H|$, then $r(C_n, H) = (n-1)(\chi(H)-1) + \sigma(H)$?
		\end{question}

	\subsection{Transversals in Latin squares}

\begin{figure}[h]
			\centering
			\begin{subfigure}{.4\textwidth}
				\centering
				\includegraphics[scale = .7]{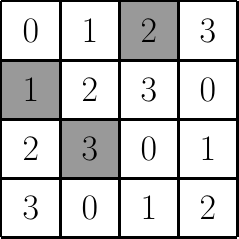}
			\end{subfigure}
			\begin{subfigure}{.4\textwidth}
				\centering
				\includegraphics[scale = .7]{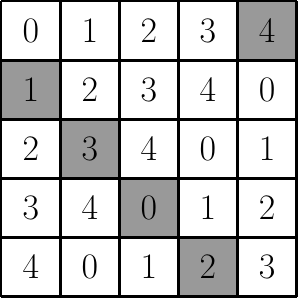}
			\end{subfigure}
			\caption{A $4 \times 4$ Latin square with a transversal of size $3$ (and with no full transversals), and a $5 \times 5$ Latin square with a full transversal}
			\label{fig:latin-square}
		\end{figure}

		A \emph{Latin square} is an $n \times n$ grid filled with $n$ symbols so that each symbol appears exactly once in each row and column. A \emph{transversal} in a Latin Square is a collection of cells that share no symbol, row, or column.
		A transversal in an $n \times n$ Latin square is said to be \emph{full} if it has size $n$. (See \Cref{fig:latin-square}.)

		The study of Latin Squares dates back to the 18$^{\text{th}}$ century, during which Euler \cite{euler1782recherches} considered Latin squares that can be decomposed into full transversals, and has been very prolific (see the surveys of Andersen \cite{andersen2013history}, Wanless \cite{wanless2011transversals}, and a forthcoming survey of Montgomery \cite{montgomery2024transversals}).
		Ryser \cite{ryser1967neuere}, Brualdi (see \cite{brualdi1991combinatorial}), and Stein \cite{stein1975transversals} made several related conjectures, which are now known in a combined form as the \emph{Ryser--Brualdi--Stein conjecture}.	
		\begin{conjecture}[Ryser--Brualdi--Stein conjecture \cite{ryser1967neuere,brualdi1991combinatorial,stein1975transversals}]
			Every $n \times n$ Latin square has a transversal of size $n-1$, and a full transversal if $n$ is odd.
		\end{conjecture}

		In a very recent breakthrough paper, Montgomery \cite{montgomery2023proof} proved the first part of this well-known conjecture, for large $n$.

		\begin{theorem}[Montgomery \cite{montgomery2023proof}] \label{thm:latin-square}
			For every large enough $n$, every $n \times n$ Latin square has a transversal of size $n-1$.
		\end{theorem}

		The substantial proof of \Cref{thm:latin-square} contains many ideas and techniques. One of them is the use of sublinear expanders, similar to $(\delta,n)$-\ssExpanders (see \Cref{def:ss-expander-var}) and to $(\lambda, d)$-\stExpanders (see \Cref{def:st-expanders}). Due to the length of the paper \cite{montgomery2023proof}, and the fact that this paper came out after the current survey was submitted for publication, we do not elaborate on the use of expanders in this proof. One can of course refer to \cite{montgomery2023proof} for more details, as well as to Montgomery's forthcoming survey \cite{montgomery2024transversals}.

	\subsection*{Acknowledgements}
		I would like to thank the anonymous referee for many useful and insightful comments. I would also like to thank Alp M\"uyesser for pointing out some typos.

\bibliography{survey}
\bibliographystyle{amsplain}

\end{document}